\documentclass[12pt]{amsart} 
\usepackage{amsmath, amsfonts, amssymb, latexsym, amscd, enumerate}
\usepackage{pictexwd,dcpic}
\usepackage{graphicx}
\usepackage{hyperref}

\newtheorem{theorem}{Theorem}[section]
\newtheorem{lemma}[theorem]{Lemma}
\newtheorem{proposition}[theorem]{Proposition}
\newtheorem{corollary}[theorem]{Corollary}
\newtheorem{remark}[theorem]{Remark}
\newtheorem{example}[theorem]{Example}
\newtheorem{definition}[theorem]{Definition}
\newtheorem{notation}[theorem]{Notation}

\numberwithin{equation}{section}

\def\zN{\mathbb N}
\def\zZ{\mathbb Z}
\def\zQ{\mathbb Q}
\def\zR{\mathbb R}
\def\zC{\mathbb C}
\def\zT{\mathbb T}
\def\cA{\mathcal A}
\def\cE{\mathcal E}
\def\cF{\mathcal F}
\def\cL{\mathcal L}
\def\cK{\mathcal K}
\def\cI{\mathcal I}
\def\cO{\mathcal O}
\def\cT{\mathcal T}

\def\fI{\mathfrak I}
\def\sq{$\text{}\hfill\square$}

\def\pf{\it Proof. \rm}
\def\st{\,\bigl\vert\,}

\newcommand{\abs}[1]{\vert#1\vert}
\newcommand{\norm}[1]{\vert\vert#1\vert\vert}
\newcommand{\innprod}[2]{\langle#1, #2\rangle}

\newcommand{\diag}[1]{\mbox{Diag}(#1)}

\title[$C^*$-algebras associated with topological group quivers I]{$C^*$-algebras associated with topological group quivers I:\\ \rm Generators, Relations and Spatial Structure}
\author{Shawn J. $\rm M^\MakeLowercase{c}$Cann}
\address{Department of Mathematics and Statistics, University of Regina, Regina, Sk, S4S 0A2}
\email{mccann1s@uregina.ca}
\begin{document}
\maketitle

\begin{abstract} Topological quivers generalize the notion of directed graphs in which
the sets of vertices and edges are locally compact (second countable) Hausdorff spaces. 
Associated to a topological quiver $Q$ is a $C^*$-correspondence, and in turn, a Cuntz-Pimsner 
algebra $C^*(Q).$ Given $\Gamma$ a locally compact group and $\alpha$ and $\beta$ 
endomorphisms on $\Gamma,$ one may construct a topological quiver $Q_{\alpha,\beta}(\Gamma)$
with vertex set $\Gamma,$ and edge set $\Omega_{\alpha,\beta}(\Gamma)=
\{(x,y)\in\Gamma\times\Gamma\st \alpha(y)=\beta(x)\}.$ In this paper, the author
examines the Cuntz-Pimsner algebra $\cO_{\alpha,\beta}(\Gamma):=C^*(Q_{\alpha,\beta}(\Gamma)).$ 
The investigative topics include a notion for topological quiver isomorphisms, generators (and their relations) 
of the $C^*$-algebras  $\cO_{\alpha,\beta}(\Gamma)$, and its spatial structure 
(i.e., colimits, tensor products and crossed products) and a few properties of its $C^*$-subalgebras.
\end{abstract}

\section{Introduction and Notation}
	\subsection{Background}The study of directed graphs have played a pivotal role in the construction and analysis of $C^*$-algebras
ultimately describing certain $C^*$-algebras in terms of easily manageable and computatable relations. That is, given
a directed graph $G=(V,E,r,s)$ with vertex set $V$, edge set $E$, range map $r$ and source map $s,$ one may 
produce a (universal) $C^*$-algebra generated by projections $\{p_v\}_{v\in V}$ and partial isometries 
$\{s_e\}_{e\in E}$ satisfying the relations\index{Graph $C^*$-algebra}
\begin{enumerate}
\item $s_e^*s_e=p_{r(e)}$ for all $e\in E$
\item $p_x=\sum_{s(e)=x} s_es_e^*$ for all $x\in X$ with $0<\abs{s^{-1}(x)}<\infty$
\item $s_es_e^*\le p_{s(e)}$ for all $e\in E.$
\end{enumerate}

Many $C^*$-algebras, such as Cuntz-algebras \cite{Cu}, Cuntz-Krieger algebras \cite{CK}, 
the Toeplitz algebra and the $n$ by $n$ matrices, are able to be defined by generators and
relations which correspond directly to such a graph.  Other $C^*$-algebras, including some that arise in representation
theory (see \cite{MRS, MRS2, PS}) and the theory of quantum spaces (see \cite{HMS, HMS2, HS}),
are of considerable importance and interest in determining how they may be represented by a graph.

Schweizer (\cite{Sz}) introduced continuous (directed) graphs and so began the analysis of spaces of vertices and edges that are topological spaces. Perhaps, the first to analyze prototypes of this structure was Deaconu 
(see \cite{Dea,Dea2,Dea3,Dea4}) who was interested in groupoid representations of Cuntz-like 
$C^*$-algebras associated to (non-invertible) maps of topological spaces, i.e. to endomorphisms of $C_0(X).$
This led Deaconu and Muhly \cite{DM} to extend the notion to branched coverings.
In order to study these $C^*$-algebras, and in particular to compute its $K$-theory, 
it was beneficial to use machinery pioneered by Pimsner in his seminal paper \cite{Pims} and realize
the $C^*$-algebra as what is now called a Cuntz-Pimsner algebra. Pimsner's
work was profoundly influenced by the theory of certain graph $C^*$-algebras (in particular, 
Cuntz-Krieger algebras) and crossed products by $\zZ$ independent of the study of ``continuous graph 
$C^*$-algebras.''

The notion of topological quivers generalizes that of directed graphs. This is 
a quintuple $Q = (X,E, r, s, \lambda),$ where $X$ and $E$ are locally compact (second countable) Hausdorff
spaces, $r$ and $s$ are continuous maps from $X$ to $E$ with $r$ open, and $\lambda = \{\lambda_x\}_{x\in E}$
is a system of Radon measures. With topological quiver, one can associate a Cuntz-Pimsner $C^*$-algebra $C^*(Q)$
that also generalizes the notion of graph $C^*$-algebra.
In a purely set-theoretic perspective, $Q$ is a directed graph and the system of measures
is to enable one to replace sums (that arise in the graph theoretic setting) with integrals. 
The term ``topological quiver'' was adopted from the ring-theoretic use of the term “quiver” 
which was introduced by Gabriel in \cite{G} to describe the graph associated to a so-called basic algebra. 

Two particularly important works were the first to investigate general structures of this nature. The first is that of
Katsura \cite{K,K1,K2,K3,K4} where the term topological graph is used to denote a structure similar to a topological quiver 
$(X,E,r,s)$\index{Topological Graph} with the addition restriction that $r$ is a local homeomorphism. If $\lambda$ is taken to be counting measures on the fibers of $r$ then it has been shown that the $C^*$-algebra in Katsura's work is isomorphic to
the $C^*$-algebra produced by considering the topological quiver $Q=(X,E,r,s,\lambda).$
The second paper is \cite{BB2} where Brenken studies topological relations $(X,E,r,s,\lambda)$ where $E$ is a closed
subset of $X\times X.$ The maps $r$ and $s$ are the projections of this set into $X.$ Such topological relations can be viewed
as multiplicity one topological quivers. \index{Topological Relations} The topological quivers considered in this paper are
examples of topological relations. 

In \cite{EaHR}, Exel, an Huef and Raeburn define $C^*$-algebras associated with a system $(B,\alpha,L)$
where $\alpha$ is an endomorphism of a unital $C^*$-algebra $B$ and $L$ is a positive linear map $L:B\to B$ such
that $L(\alpha(a)b) = aL(b)$ for all $a, b\in B$ called a \emph{transfer operator}. In fact, the $C^*$-algebra they generate is a Cuntz-Pimsner algebra
and under certain restrictions, a $C^*$-algebra associated with a topological quiver;\index{Transfer Operator}
in particular, when $B=C(\zT^d)$ the continuous function on the $d$-torus,
$F\in M_d(\zZ)$ and $\alpha$ is the endomorphism 
$$\alpha(f)(e^{2\pi it})=f(e^{2\pi i Ft})$$
for $f\in C(\zT^d)$ and $t\in\zR^d.$ 
Furthermore, Yamashita \cite{Yam} considers topological relations
of the form $Q_{n,m}(\zT)=(\zT,\Omega_{n,m}(\zT),r,s,\lambda)$ where $n,m\in\zZ,$ $\gcd(n,m)=1$ and
$$\Omega_{n,m}(\zT)=\{(x,y)\in\zT\times\zT\st y^n=x^m\}.$$
He then presents the $C^*$-algebra associated with $Q_{n,m}(\zT)$
as a universal $C^*$-algebra with certain generators and relation; that is, the unversal $C^*$-algebra generated by
a unitary $U$ and isometries $S_1,...,S_n$ such that
$$US_i=S_{i+1}\,(1\le i\le n-1),\qquad US_n=S_1U^m,\qquad \sum_{i=1}^n S_iS_i^*=1.$$
We generalize the notions considered in \cite{EaHR} and \cite{Yam} and provide an extensive survey.

In Section 2, the necessary background material is provided in order to define topological quivers and the $C^*$-algebra
associated with a topological quiver. We then develop a method for determining a presentation of the universal $C^*$-algebra
and define a notion for a topological quiver isomorphism and prove this morphism gives an isomorphism between the associated
$C^*$-algebras. Section 3 then contains the definition of a particular topological quiver, called
a topological group relation, defined using a locally compact group $\Gamma$ and endomorphisms $\alpha$ and $\beta;$
that is, 
$$Q_{\alpha,\beta}(\Gamma)=(\Gamma, \Omega_{\alpha,\beta}(\Gamma), r,s,\lambda)$$
where 
$$\Omega_{\alpha,\beta}(\zT^d)=\{(x,y)\in\Gamma\times\Gamma\st\alpha(y)=\beta(x)\}.$$
For instance, if $F,G\in M_d(\zZ)$ where $\det F, \det G\ne0,$ then
$$Q_{F,G}(\zT)^d=(\zT^d,\Omega_{F,G}(\zT^d),r,s,\lambda)$$
where 
$$\Omega_{F,G}(\zT^d)=\{(x,y)\in\zT^d\times\zT^d\st \sigma_F(y)=\sigma_G(x)\}$$
and $\sigma_F(e^{2\pi it})=e^{2\pi iFt}$ for $t\in\zR^d.$ 
These topological group relations are a generalization of the topological quivers studied in \cite{EaHR} and \cite{Yam} and
they also share many similar properties with the Cuntz algebra. We also consider many other examples of topological
group relations and also presents the $C^*$-algebras, $\cO_{F,G}(\zT^d)$ 
associated with $Q_{F,G}(\zT^d)$ as a universal $C^*$-algebra
generated by commuting unitaries $\{U_j\}_{j=1}^d$ and an isometry $S$ satisfying certain relations:
\begin{enumerate}
\item $S^*U^\nu S=\delta_{\nu}^{\nu^\prime},$
\item $U_j^{a_j}S=SU^{G_j},$ for all $j=1,...,d$ and
\item $1=\sum_{\nu\in\fI(F)} U^\nu SS^*U^{-\nu}$
\end{enumerate}
where $U^\nu$ denotes $\prod_{j=1}^dU_j^{\nu_j}$ and $\fI(F)=\{\nu=(\nu_j)\in\zZ\st 0\le \nu_j\le a_j-1\}.$ 
In Section 4, we investigate a few properties of 
$C^*$-subalgebras of $\cO_{F,G}(\zT^d)$ while in Section 5, we examine
the spatial structure of $\cO_{F,G}(\zT^d)$ and present $\cO_{F,G}(\zT^d)$ as the crossed product by an endomorphism
of a colimit. In \cite{S}, Stacey defines the \emph{crossed product by an endomorphism},\index{Crossed Product by an Endomorphism} $\rho,$ 
as the universal $C^*$-algebra $B\rtimes_\rho\zN$ generated by (a copy of) $B$ and an isometry $s$ such that $\rho(b)=sbs^*$ for all $b\in B.$ Ultimately, it is shown that 
$$\cO_{F,G}(\zT^d)=(\mbox{colim}_k M_{N^k}(C(\zT^d)))\rtimes_\rho\zN$$ where the isometry $s$ is, none other than, $S.$

	\subsection{Notation}The sets of natural numbers, integers, rationals numbers, real numbers and complex numbers will be denoted by
$\zN$, $\zZ$, $\zQ$, $\zR$, and $\zC,$ respectively. Convention: $\zN$ does not contain zero.
Finally, $\zZ_p$ denotes the abelian group $\zZ/p\zZ=\{0,1,...,p-1\mod p\}$ and 
$\zT$ denotes the torus $\{z\in\zC\st \abs{z}=1\}.$ Whenever convenient, view $\zZ_p\subset\zT$ by 
$\zZ_p\cong\{z\in\zT\st z^p=1\}.$

For a topological space $Y$, the closure of $Y$ is denoted $\overline{Y}.$ Given a locally compact Hausdorff space $X$, let
$C(X)$ be the continuous complex functions on $X,$
$C_0(X)$ be the continuous complex functions on $X$ vanishing at infinity,
and $C_c(X)$ be the continuous complex functions on $X$ with compact support.
The supremum norm, denoted $\norm{\cdot}_\infty,$ is defined by
$$\norm{f}_\infty=\sup_{x\in X}\{\abs{f(x)}\}$$
for each continuous map $f:X\to\zC.$ For a continuous function $f\in C_c(X),$ denote the open support of $f$ by 
$\mbox{osupp }f=\{x\in X\st f(x)\ne 0\}$ and the support of $f$ by $\mbox{supp }f=\overline{\mbox{osupp} f}.$  

For $C^*$-algebras $A$ and $B$, $A$ is isomorphic to $B$ will be written $A\cong B;$ for example, we use
$C(\zT^d)\otimes M_{N}(\zC)\cong M_{N}(C(\zT^d)).$ 
Moreover, $A^{\oplus n}$ denotes the $n$-fold direct sum $A\oplus\cdots\oplus A.$ 
Given a group $\Gamma$ and a ring $R$, a normal subgroup, $N$, of $\Gamma$ is denoted $N\lhd\Gamma$ and
an ideal, $I$, of $R$ is denoted $I\lhd R.$ Note if $R$ is a $C^*$-algebra then the term ideal denotes a closed two-sided ideal. 
Furthermore, $\mbox{End}(\Gamma)$ ($\mbox{End}(R)$) and $\mbox{Aut}(\Gamma)$
($\mbox{Aut}(R)$) denotes the set of endomorphisms of $\Gamma$ ($R$) and automorphisms of $\Gamma$ $(R$), respectively. 
For a map $\gamma:\Gamma\to\mbox{Aut}(A),$ the \emph{fixed point set}\index{Fixed Point Set} 
is denoted $A^\gamma$ and defined by 
$$A^\gamma=\{a\in A\st \gamma(g)(a)=a\mbox{ for each }g\in \Gamma\}.$$

Let $\alpha\in C(X)$ then $\alpha^\#\in\mbox{End(C(X))}$ denotes the endomorphism of $C(X)$ defined by
$$\alpha^\#(f)=f\circ\alpha\qquad\mbox{for each $f\in C(X)$}.$$
Let $S$ be a set and define the Kronecker delta function $\delta:S\times S\to\{0,1\}$ by
$$\delta_s^r:=\delta(s,r)=\begin{cases}
0&\mbox{if }s\ne r\\
1&\mbox{if }s=r
\end{cases}.$$
\section{Preliminairies}\subsection{Hilbert $C^*$-modules}\begin{definition}\label{Hbmod}\cite{lan} \rm  If $A$ is a $C^*$-algebra, then a \emph{(right) Hilbert $A$-module}\index{Hilbert $C^*$-module} is a Banach space $\cE_A$
together with a right action of $A$ on $\cE_A$ and an $A$-valued inner product $\innprod{\cdot}{\cdot}_A$
satisfying
\begin{enumerate}
\item $\innprod{\xi}{\eta a}_A=\innprod{\xi}{\eta}_A a$
\item $\innprod{\xi}{\eta}_A =\innprod{\eta}{\xi}_A^*$
\item $\innprod{\xi}{\xi}\ge 0$ and $\norm{\xi}=\norm{\innprod{\xi}{\xi}_A^{1/2}}_A$
\end{enumerate}
for all $\xi$, $\eta\in\cE_A$ and $a\in A$ (if the context is clear, we denote $\cE_A$ simply by $\cE$). 
For Hilbert $A$-modules $\cE$ and $\cF$, call a function $T:\cE\to\cF$ \emph{adjointable}
\index{Hilbert $C^*$-module! Adjointable Operator, $\cL(\cE,\cF)$}
\index{Adjointable Operator, $\cL(\cE,\cF)$}
if there is a function $T^*:\cF\to\cE$ such that
$\innprod{T(\xi)}{\eta}_A=\innprod{\xi}{T^*(\eta)}_A$ for all $\xi\in\cE$ and $\eta\in\cF$.
Let $\cL(\cE,\cF)$ denote the set of adjointable ($A$-linear) operators from $\cE$ to $\cF$. 
If $\cE=\cF$, then $\cL(\cE):=\cL(\cE,\cE)$ is a $C^*$-algebra.
Let $\cK(\cE, \cF)$ denote the closed two-sided ideal of \emph{compact operators}
\index{Hilbert $C^*$-module! Compact Operators, $\cK(\cE,\cF)$} 
\index{Compact Operators, $\cK(\cE,\cF)$}
given by
$$\cK(\cE,\cF):=\overline{\mbox{span}}\{\theta_{\xi,\eta}^{\cE,\cF}\st\xi\in\cE,\,\eta\in\cF\}$$
where $$\theta_{\xi,\eta}^{\cE,\cF}(\zeta)=\xi\innprod{\eta}{\zeta}_A\qquad\mbox{for each $\zeta\in\cE$}.$$ 
Similarly, $\cK(\cE):=\cK(\cE,\cE)$ and $\theta_{\xi,\eta}^\cE$ (or $\theta_{\xi,\eta}$ if understood) denotes 
$\theta_{\xi,\eta}^{\cE,\cE}$.
For Hilbert $A$-module $\cE$, the linear span of $\{\innprod{\xi}{\eta}\st\xi,\eta\in\cE\}$, denoted $\innprod{\cE}{\cE}$,
once closed is a two-sided ideal of $A$. Note that $\cE\innprod{\cE}{\cE}$ is dense in $\cE.$ 
The Hilbert module $\cE$ is called \emph{full}\index{Hilbert $C^*$-module! Full} if $\innprod{\cE}{\cE}$ is dense in $A$. 
The Hilbert module $A_A$ refers to the Hilbert module $A$ over itself, where $\innprod{a}{b}=a^*b$ for all $a,b\in A$.

An \emph{algebraic generating set}\index{Hilbert $C^*$-module! Algebraic Generating Set} for $\cE$ is a subset $\{u_i\}_{i\in\cI}\subset \cE$ for some indexing set $\cI$ such that $\cE$ equals the linear span of $\{u_i\cdot a\st\ i\in\cI, a\in A\}.$ 
\end{definition}

\begin{definition}\label{ONB}\cite{KW} \rm A subset $\{u_i\}_{i\in\cI}\subset \cE$ is called a \emph{basis}\index{Basis}
provided the following reconstruction formula holds for all $\xi\in\cE:$
$$ \xi=\sum_{i\in\cI} u_i\cdot\innprod{u_i}{\xi}\qquad(\mbox{in }\cE,\norm{\cdot}.)$$
If $\innprod{u_i}{u_j}=\delta_i^j$ as well, call $\{u_i\}_{i\in\cI}$ an \emph{orthonormal basis}
\index{Basis! Orthonormal} of $\cE$.
\end{definition}  

\begin{remark}\rm The preceding definition is in accordance with the finite version in \cite{KW}, 
but many other versions exist such as in
\cite{EaHR} where $\{u_i\}_{i=1}^n$ is called a finite Parseval frame, or in \cite{Yam} where this is taken as the definition
for \it finitely generated.\rm\index{Hilbert $C^*$-module! Finitely Generated} 
There has been substantial work done on similar frames (see \cite{HJLM}).
\end{remark}

\begin{definition}\label{Hbmorph}\cite{BB4, BB5} \rm If $A$ and $B$ are $C^*$-algebras, then an \emph{$A-B$ $C^*$-correspondence}\index{C$\mbox{}^*$-correspondence} $\cE$  is a right Hilbert 
$B$-module $\cE_B$ together with a left action of $A$ on $\cE$ given by a $*$-homomorphism $\phi_A:A\to\cL(\cE)$, 
$a\cdot\xi=\phi_A(a)\xi$ for $a\in A$ and $\xi\in\cE$. We may occasionally write, $_A\cE_B$ to denote an $A-B$ 
$C^*$-correspondence and $\phi$ instead of $\phi_A$.
Furthermore, if $_{A_1}\cE_{B_1}$ and $_{A_2}\cF_{B_2}$ are $C^*$-correspondences, then
a \emph{morphism}\index{C$\mbox{}^*$-correspondence! Morphism} 
$(\pi_1, T, \pi_2):\cE\to\cF$ consists of $*$-homomorphisms $\pi_i:A_i\to B_i$ and a linear map
$T:\cE\to\cF$ satisfying
\begin{enumerate}
\item[(i)] $\pi_2(\innprod{\xi}{\eta}_{A_2})=\innprod{T(\xi)}{T(\eta)}_{B_2}$
\item[(ii)] $T(\phi_{A_1}(a_1)\xi)=\phi_{B_1}(\pi_1(a_1))T(\xi)$
\item[(iii)] $T(\xi)\pi_2(a_2)=T(\xi a_2)$
\end{enumerate} for all $\xi,\eta\in\cE$ and $a_i\in A_i$.
\end{definition}

\begin{notation}\rm When $A=B$, we refer to $_A\cE_A$ as a $C^*$-correspondence over $A$. For $\cE$ a 
$C^*$-correspondence over $A$ and $\cF$ a $C^*$-correspondence over $B$, a morphism 
 $(\pi,T,\pi):\cE\to\cF$ will be denoted by $(T,\pi)$. 
\end{notation}

\begin{definition}\cite{MT} \rm If $\cF$ is the Hilbert module $ _CC_C$ where $C$ is a $C^*$-algebra with the inner product 
$\innprod{x}{y}_B=x^*y$ then call a morphism $(T,\pi):$ $_A\cE_B\to C$ of Hilbert
modules a \emph{representation}\index{C$\mbox{}^*$-correspondence! Representation} of $_A\cE_B$ into $C.$ 
\end{definition}

\begin{remark}\rm Note that a representation of $_A\cE_B$ need only satisfying $(i)$ and $(ii)$ 
of definition \ref{Hbmorph} as it was unnecessary to require (iii) since 
\begin{align*}
\norm{T(\xi)\pi(a)-T(\xi a)}^2&=\norm{(T(\xi)\pi(a)-T(\xi a))^*(T(\xi)\pi(a)-T(\xi a))}\\
&=\norm{\pi(a)^*T(\xi)^*T(\xi)\pi(a)-T(\xi a)^*T(\xi)\pi(a)\\
&\qquad-\pi(a)^*T(\xi)^*T(\xi a) + T(\xi a)^*T(\xi a)}\\
&=\norm{\pi(a^*\innprod{\xi}{\xi}a-\innprod{\xi a}{\xi}a-a^*\innprod{\xi}{\xi a}+\innprod{\xi a}{\xi a})}\\
&=\norm{\pi(a^*\innprod{\xi}{\xi}a-a^*\innprod{\xi}{\xi}a-a^*\innprod{\xi}{\xi}a+a^*\innprod{\xi}{\xi}a)}=0
\end{align*} by condition (i).
\end{remark}

A morphism of Hilbert modules $(T,\pi):\cE\to\cF$ yields a $*$-homomorphism $\Psi_T:\cK(\cE)\to\cK(\cF)$ by
$$\Psi_T(\theta_{\xi,\eta}^\cE)=\theta_{T(\xi),T(\eta)}^\cF$$
for $\xi,\eta\in\cE$ and if $(S,\sigma):\mathcal D\to\cE$, and $(T,\pi):\cE\to\cF$ are morphisms of Hilbert modules then
$\Psi_T\circ\Psi_S=\Psi_{T\circ S}$. In the case where $\cF=B$ a $C^*$-algebra, we may first identify $\cK(B)$ as $B$,
and a representation $(T,\pi)$ of $\cE$ in a $C^*$-algebra $B$ yields a $*$-homomorphism $\Psi_T:\cK(\cE)\to B$ given
by $$\Psi_T(\theta_{\xi,\eta})=T(\xi)T(\eta)^*.$$

\begin{definition}\cite{MT} \rm For a $C^*$-correspondence $\cE$ over $A$, denote the ideal $\phi^{-1}(\cK(\cE))$ of $A$ by $J(\cE),$\index{J$(\cE)$} and let $J_\cE=J(\cE)\cap(\ker \phi)^\perp$\index{J$\mbox{}_\cE$} where $(\ker\phi)^\perp $ is the 
ideal $\{a\in A\st ab=0 \mbox{ for all }b\in\ker\phi\}$ .
If $_A\cE_A$ and $_B\cF_B$ are $C^*$-correspondences over $A$ and $B$ respectively and $K\lhd J(\cE)$, a morphism
$(T,\pi):\cE\to\cF$ is called \emph{coisometric on $K$}\index{C$\mbox{}^*$-correspondence! Representation! Coisometric on $K$} if $$\Psi_T(\phi_A(a))=\phi_B(\pi(a))$$
for all $a\in K$, or just \emph{coisometric}, if $K=J(\cE)$.
\end{definition}

\begin{notation}\rm We denote $C^*(T,\pi)$ to be the $C^*$-algebra generated by $T(\cE)$ and $\pi(A)$ where
$(T,\pi):\cE\to B$ is a representation of $_A\cE_A$ in a $C^*$-algebra $B$. Furthermore, if $\rho:B\to C$ is a 
$*$-homomorphism of $C^*$-algebras, then $\rho\circ (T,\pi)$ denotes the representation $(\rho\circ T,\rho\circ\pi)$ of $\cE$.
\end{notation}

\begin{definition}\label{CPDefs}\cite{MT} \rm 
A morphism $(T_\cE,\pi_\cE)$ coisometric on an ideal $K$ is said to be \emph{universal}
if whenever $(T,\pi):\cE\to B$ is a representation coisometric on $K$, there exists a $*$-homomorphism 
$\rho:C^*(T_\cE,\pi_\cE)\to B$ with $(T,\pi)=\rho\circ(T_\cE,\pi_\cE)$.  The universal $C^*$-algebra 
$C^*(T_\cE,\pi_\cE)$ is called the \emph{relative Cuntz-Pimsner algebra} 
of $\cE$ determined by the ideal $K$ and denoted by $\cO(K,\cE)$. If $K=0$, then $\cO(K,\cE)$ is denoted
by $\cT(\cE)$ and called the \emph{universal Toeplitz $C^*$-algebra} for $\cE$. We denote $\cO(J_\cE,\cE)$ by $\cO_\cE$. 
\end{definition}

In \cite{K}, Katsura provides a sufficient condition to guarantee the nuclearity of the Cuntz-Pimsner algebra $\cO_{\cE}$: 

\begin{proposition}\label{KatNuc}\cite[Corollary 7.4]{K} \rm  
Let $\cE$ be a $C^*$-correspondence over $C^*$-algebra $A.$
If $A$ is nuclear, then $\cO_{\cE}$ is nuclear.
\end{proposition}

	\subsection{Topological Quivers}

Many $C^*$-algebras are achieved as graph $C^*$-algebras \cite{Kum, McThesis} such as $0$, $\zC^d$ ($d\in\zN$),
$M_n(\zC)$ ($n\in\zN$), $\cK$ (the compact operators on a separable Hilbert space), the Cuntz algebras $\cO_n$ ($n\in\zN$)
and the Toeplitz Algebra $\cT.$ It was proven in \cite{FLR} for sinkless graph $G$ and generalized in \cite{BB2.5} 
that $C^*(G)$ is a (relative) Cuntz-Pimsner algebra.

A certain notion of ``continuous graph $C^*$-algebra'' exists when $X$ and $E$ are infinite (perhaps nondiscrete) sets, yet a few extra manageable properties for $X$, and $E$ are needed. 
This notion is what is called a topological quiver which may be found in \cite{BB3, BB4, BB5, MT}.

\begin{definition}\label{TQ}\cite{MT} \rm A \emph{topological quiver}\index{Topological Quiver} 
(or \emph{topological directed graph}\index{Topological Directed Graph}) $Q=(X,E,Y,r,s,\lambda)$ is a diagram
$$\begindc{\commdiag}[5]
\obj(10,0){$E$}
\obj(0,0){$X$}
\obj(20,0){$Y$}
\mor{$E$}{$X$}{$s$}[-1,0]
\mor{$E$}{$Y$}{$r$}
\enddc$$
where $X,E,$ and $Y$ are second countable locally compact Hausdorff spaces, $r$ and $s$ are continuous maps with $r$ open, along with a family $\lambda=\{\lambda_y\vert y\in Y\}$ of Radon measures on $E$ satisfying
\begin{enumerate}
\item $\mbox{supp }\lambda_y=r^{-1}(y)$ for all $y\in Y$, and
\item $y\mapsto\lambda_y(f)=\int_Ef(\alpha)d\lambda_y(\alpha)\in C_c(Y)$ for $f\in C_c(E).$
\end{enumerate}
If $X=Y$ then write $Q=(X,E,r,s,\lambda)$ in lieu of $(X,E,X,r,s,\lambda).$
\end{definition}

\begin{remark}\rm 
It is important to remark that in the literature pertaining to graph $C^*$-algebras and their generalizations, 
various authors interchange the roles of the maps $r$ and $s$ in their definitions. 
Definition \ref{TQ} agrees with that used in most papers on graph $C^*$-algebras (cf. \cite{BHRS, BPRS, BB2, BB3, BB4, KPR, KPRR, MT}). However, Definition \ref{TQ} differs from that used in the higher rank graph algebras of Kumjian and 
Pask \cite{KP1, KP2} and from that used in the topological graph algebras of Katsura \cite{K1, K2, K3, K4} where Katsura uses the term topological graph to denote a more restrictive concept. 
\end{remark}

Given a topological quiver $Q=(X,E,Y,r,s,\lambda)$, one may associate a correspondence  
$\cE_Q$ of the $C^*$-algebra $C_0(X)$ to the $C^*$-algebra $C_0(Y)$. Define left and right actions 
$$(a\cdot\xi\cdot b)(e)=a(s(e))\xi(e)b(r(e))$$
by $C_0(X)$ and $C_0(Y)$ respectively on $C_c(E)$. Furthermore, define the $C_c(Y)$-valued inner product
$$\innprod{\xi}{\eta}(y)=\int_{r^{-1}(y)}\overline{\xi(\alpha)}\eta(\alpha)d\lambda_y(\alpha)$$ 
for $\xi,\eta\in C_c(E)$, $y\in Y,$ and let $\cE_Q$\index{C$\mbox{}^*$-correspondence! Associated with a Topological Quiver} be the completion of $C_c(E)$ with respect to the norm 
$$\norm{\xi}=\norm{\innprod{\xi}{\xi}^{1/2}}_\infty=\norm{\lambda_y(\abs{\xi}^2)}_\infty^{1/2}.$$ 

\begin{definition}\rm Given topological quiver $Q$ over a space $X$, define the $C^*$-algebra, $C^*(Q)$
\index{C$\mbox{}^*$-algebra Associated with a Topological Quiver}\index{C$\mbox{}^*(Q)$}
\index{Topological Quiver! $C^*$-algebra Associated with} 
associated with $Q$ to be the Cuntz-Pimnser $C^*$-algebra $\cO_{\cE_Q}$\index{Cuntz-Pimsner Algebra} 
of the correspondence $\cE_Q$ over $A=C_0(X)$.
\end{definition}

\begin{proposition}\cite[Proposition 2.21]{MS0}\label{InjA} \rm 
If the left action, $\phi,$ is injective, then the universal map $i_A : A\to \cO_{\cE_Q}$ is injective.
\end{proposition}

\begin{proposition}\label{Nuc}\rm $C^*(Q)$ is nuclear for any topological quiver $Q.$\\
\pf Since abelian $C^*$-algebras are nuclear, so is $A=C_0(X)$ and hence by Proposition \ref{KatNuc},
$C^*(Q)$ is nuclear.\\\sq 
\end{proposition}

\begin{example}\rm \mbox{} 
\begin{enumerate} 
\item Let $G=(X,E,r,s)$ be a directed graph and let $\lambda_y$ be counting measure on $r^{-1}(y)$. 
Then $G$ along with this family of Radon measures $\lambda$ on $E$ becomes a topological quiver $Q=(X,E,r,s,\lambda)$.  
As we have already mentioned, $C^*(G)\cong C^*(Q).$
\item In \cite{K1}, Katsura uses topological graph $G=(X,E,r,d)$ to be the more restrictive situation 
using two continuous maps $r,d:E\to X$ with $d$ a local homeomorphism. We may set $\lambda_x$ to be counting measure on $d^{-1}(x)$ since $d^{-1}(x)$ must be a discrete space.
Then $Q=(X,E,d,r,\lambda)$ is a topological quiver. Furthermore, Katsura defines an associated $C^*$-algebra $C^*(G)$ which is
isomorphic to $C^*(Q)$ \cite[Definition 2.9 and Definition 2.10]{K1}.\index{Topological Graph}
\item In \cite{BB2}, Brenken defines a $C^*$-algebra associated to a closed relation, $\alpha\subset X\times X,$ on a topological space $X$ with $\pi_i$ the usual projection map onto the $i$-th coordinate. There is a Radon measure $\mu_x$ such that $\mbox{supp }\mu_x=\pi_2^{-1}(x)$ for all $x\in X$. The Cuntz-Pimsner algebra of the associated correspondence, $C^*(\alpha)$ 
is isomorphic to $C^*(Q)$ where $Q$ is the topological quiver $(X,\alpha, \pi_2, \pi_1, \mu).$ Such topological quivers are
referred to as \emph{topological relations}.\index{Topological Relations}
\end{enumerate}
\end{example}

\begin{theorem}\label{BasisGen}\rm Let $Q=(X, E, r,s,\lambda)$ be a topological quiver with $X$ compact.
If $\cE_Q$ has a finite orthonormal basis\index{Basis! Orthonormal} $\{u_i\}_{i=1}^n\subset C_c(E)$ and $\phi$ is injective then 
$C^*(Q)$ is the universal $C^*$-algebra generated by $C(X)$ and elements $\{S_\xi\}_{\xi\in\fI},$ where $\fI$ is 
the module generated by $\{u_i\}_{i=1}^n,$ satisfying the relations:
\begin{enumerate}
\item $S_{\alpha\xi+\beta\eta}=\alpha S_\xi+\beta S_\eta,$
\item $a\cdot S_\xi\cdot b=S_{a\cdot\xi\cdot b},$
\item$S_{u_i}^*S_{u_j}=\delta_i^j1\in C(X)$, and
\item $\sum_{i=1}^n S_{u_i}S_{u_i}^*=1$
\end{enumerate}
for all $\alpha,\beta\in\zC$, $a,b\in C(X)$ and $\xi,\eta\in\fI$. Moreover, the Toeplitz-Pimsner algebra 
\index{Toeplitz-Pimsner Algebra}$\cT(Q):=\cT(\cE_Q)$ is the 
universal $C^*$-algebra satisfying (1)-(3).\\
\pf Let $\cA$ be the universal $C^*$-algebra generated by $C(X)$ and $\{S_\xi\st\xi\in\fI\}$ satisfying the relations:
\begin{enumerate}
\item $S_{\alpha\xi+\beta\eta}=\alpha S_\xi+\beta S_\eta,$
\item $a\cdot S_\xi\cdot b=S_{a\cdot\xi\cdot b},$
\item$S_{u_i}^*S_{u_j}=\delta_i^j1$, and
\item $\sum_{i=1}^n S_{u_i}S_{u_i}^*=1$
\end{enumerate}
for all $\alpha,\beta\in\zC$, $a,b\in C(X)$ and $\xi,\eta\in\fI.$ To show the existence of $\cA$, calculate
\begin{align*}
\norm{S_\xi}&=\norm{S_{\sum_{i=1}^nu_i\innprod{u_i}{\xi}}}\\
&=\norm{\sum_{i=1}^nS_{u_i}\innprod{u_i}{\xi}}\\
&\le\sum_{i=1}^n\norm{S_{u_i}}\norm{\innprod{u_i}{\xi}}_\infty
\end{align*} 
for each $\xi\in\fI$ and furthermore, $\norm{S_{u_i}}^2=\norm{S_{u_i}^*S_{u_i}}=1.$ Since 
$\norm{\innprod{u_i}{\xi}}_\infty$ is bounded in $C(X),$ $\cA$ exists as a universal $C^*$-algebra.

Now if $a\in C(X)$ then 
$$\sum_{i=1}^n \theta_{a\cdot u_i,u_i}(\xi)(e)=\sum_{i=1}^n (a\cdot u_i)\innprod{u_i}{\xi}(e)=a(s(e))\xi(e)=(\phi(a)\xi)(e);$$
that is, $J(\cE_Q)=C(X).$
Then let $(T_u,\pi_u)$ be the universal representation of $\cE_Q$ in $C^*(Q)$ coisometric
on $J_{\cE_Q}=J(\cE_Q)$ (since $\phi$ is injective.) Setting $S_\xi=T_u(\xi)$, we have
$S_{\alpha\xi+\beta\eta}=\alpha S_\xi+\beta S_\eta,$ $a\cdot S_\xi\cdot b=S_{a\cdot\xi\cdot b},$ and
$S_{u_i}^*S_{u_j}=\innprod{u_i}{u_j}=\delta_i^j1$
for all $\alpha,\beta\in\zC$, $a,b\in C(X)$ and $\xi,\eta\in\fI.$ 
Finally, the reconstruction formula can easily be rewritten as
$$\sum_{i=1}^n\theta_{u_i,u_i}=1\in C(X)=J_{\cE_Q};$$
hence, (4) is evident. Thus, there exists a homomorphism $\rho:\cA\to C^*(Q)$ such that 
$$S_\xi\mapsto T_u(\xi)$$
and $$a\mapsto \pi_u(a)$$
for each $\xi\in\fI$ and $a\in C(X).$

Conversely, given $\{S_\xi\}_{\xi\in\fI}$ satisfying (1)-(4), we construct a coisometric representation $(T,\pi)$ of $\cE_Q$ into $\cA$; that is, let $T(\xi)=S_{\xi}$ for each $\xi\in\cE_Q (=\fI$, since $\{u_i\}_{i=1}^n$ is a basis)
and let $\pi$ be the inclusion of $C(X)$ into the universal $C^*$-algebra stated above. 
With this in mind, (1) and (2) imply that $T$ is a linear map that satisfies
$$T(a\cdot\xi\cdot b)=a\cdot T(\xi)\cdot b.$$
Furthermore, if $\xi,\eta\in\cE_Q$ then
\begin{align*}
T(\xi)^*T(\eta)&=T(\sum_{i=1}^nu_i\innprod{u_i}{\xi})^*T(\sum_{j=1}^nu_j\innprod{u_j}{\eta})\\
&=[\sum_{i=1}^n\innprod{\xi}{u_i}S_{u_i}^*][\sum_{j=1}^nS_{u_j}\innprod{u_j}{\eta}]\\
&=\sum_{i=1}^n\innprod{\xi}{u_i}\innprod{u_i}{\eta}=\sum_{i=1}^n\innprod{u_i\innprod{u_i}{\xi}}{\eta}\\
&=\innprod{\sum_{i=1}^nu_i\innprod{u_i}{\xi}}{\eta}=\innprod{\xi}{\eta}.
\end{align*}
Finally, if $a\in C(X)$ then 
$$\Psi_T(\phi(a))=\Psi_T(\sum_{i=1}^n\theta_{a\cdot u_i,u_i})=\sum_{i=1}^n S_{a\cdot u_i}S_{u_i}^*=a\sum_{i=1}^n S_iS_i^*=a.$$
Thus, $\Psi_T(\phi(a))=a$ and hence, $(T,\pi)$ is coisometric. Thus, there exist a homomorphism $\tau:C^*(Q)\to\cA$
such that $\tau\circ(T_u,\pi_u)\to (T,\pi)$ where $(T_u,\pi_u)$ is the universal coisometric representation of $\cE_Q$ 
in $C^*(Q)$. Thus,
$$\tau\circ\rho(S_\xi)=\tau(T_u(\xi))=T(\xi)=S_\xi$$
and $$\tau\circ\rho(a)=\tau(\pi_u(a))=\pi(a)=a$$
for each $\xi\in\cE_Q$ and $a\in A$; that is, $\rho$ and $\tau$ are isomorphisms.\\\sq
\end{theorem}

\begin{remark}\rm The existence of a basis $\{u_i\}_{i=1}^n$ forces $\sum_{i=1}^n\theta_{u_i,u_i}=1$
and hence, $X$ must be compact and $\phi$ must be injective.
\end{remark}





	\subsection{Topological Quiver Isomorphisms}A notion for a topological quiver isomorphism based on the definition of a directed graph isomorphism is introduced.

\begin{definition}\rm  Let $Q=(X,E,Y,r,s,\lambda)$ and $Q^\prime=(X^\prime,E^\prime,Y^\prime,r^\prime,s^\prime,\lambda^\prime)$ be topological quivers.\index{Topological Quiver}  If there exist homeomorphisms $\varphi_1:X\to X^\prime$, $\varphi_2:Y\to Y^\prime$, and $\psi:E\to E^\prime$ where $\psi$ is measurable and there exists a family of constants $\{h_y>0\st y\in Y\}$ such that
\begin{enumerate}
\item $\varphi_1(s(e))=s^\prime(\psi(e))$
\item $\varphi_2(r(e))=r^\prime(\psi(e))$ 
\item $\lambda^\prime_{\varphi(y)}(B)=\int_B h_y\,d(\lambda_y\circ\psi^{-1})$ 
\end{enumerate}
for each $e\in E,$ measurable $B\subseteq E$ and $y\in Y,$ then say $Q$ and $Q^\prime$ are \emph{isomorphic}
\index{Topological Quiver! Isomorphism} and write $Q\cong Q^\prime$.   
\end{definition}

\begin{remark}\rm Recall that a measurable map $\psi:E\to E^\prime$ is one that satisfies $\psi^{-1}(B^\prime)$ is measurable
in $E$ for each measurable subset $B^\prime\subseteq E^\prime$ (see \cite{W}.) For any homeomorphism $\psi:E\to E^\prime$ and compact $K\subset E^\prime,$
assume $\psi^{-1}(K)$ is covered by the open collection $\{U_i\}_{i\in\cI}.$ Then $\psi^{-1}(K)\subset\cup_{i\in\cI} U_i$
implies $K\subset\cup_{i\in\cI} \psi(U_i).$ Hence, $K$ is covered by the open sets $\psi(U_i)$ and so, $K\subset\cup_{j=1}^n\psi(U_{i_j})$ for some open subcover $\{\psi(U_{i_j})\}_{j=1}^n.$ Finally,
$\psi^{-1}(K)\subset\cup_{j=1}^n U_{i_j};$ that is, $\psi^{-1}(K)$ is compact and thus, $\psi$ is a proper map. 
Note assuming $\psi$ is open, injective and continuous is enough to show that $\psi$ is a proper map. 
\end{remark} 

\begin{theorem}\rm \label{QuiverIso} Let $Q=(X,E,r,s,\lambda)$ and $Q^\prime=(X^\prime,E^\prime,r^\prime,s^\prime,\lambda^\prime)$ be topological quivers. If $Q\cong Q^\prime$ by $(\varphi,\psi, h)$, then $\cO(K,\cE_Q)\cong \cO(K^\prime,\cE_{Q^\prime})$ for any closed ideal $K^\prime$ of $J(\cE_Q)$ where $K=\varphi^\#(K^\prime)$; that is, given a homeomorphism 
$\varphi:X\to X^\prime,$ a measurable homeomorphism $\psi:E\to E^\prime$ and a family of constants $\{h_y>0\st y\in X\}$  
satisfying 
\begin{enumerate}
\item $\varphi_1(s(e))=s^\prime(\psi(e))$
\item $\varphi_2(r(e))=r^\prime(\psi(e))$ 
\item $\lambda^\prime_{\varphi(y)}(B)=\int_B h_y\,d(\lambda_y\circ\psi^{-1})$ 
\end{enumerate}
for each $e\in E,$ measurable $B\subseteq E$ and $y\in Y$, then $\cO(K,\cE_Q)\cong \cO(K^\prime,\cE_{Q^\prime})$.\\
\pf Define $T_0:C_c(E^\prime)\to C_c(E)$ and $\pi:C_0(X^\prime)\to C_0(X)$ by
$$T_0(\xi^\prime)(e)=h_{r(e)}^{1/2}\xi^\prime(\psi(e))\qquad\mbox{and}\qquad\pi(a^\prime)=a^\prime\circ\varphi$$
for all $\xi^\prime\in C_c(E^\prime)$ and $a^\prime\in C_0(X^\prime)$. Since $\psi$ is proper, $T_0(\xi^\prime)\in C_c(E)$
for each $\xi^\prime\in C_c(E^\prime).$ For $\xi^\prime,\eta^\prime\in C_c(E^\prime),$ $\alpha,\beta\in\zC$ and $e\in E,$
\begin{align*}
T_0(\alpha\xi^\prime+\beta\eta^\prime)(e)&=h_{r(e)}^{1/2}(\alpha\xi^\prime+\beta\eta^\prime)(\psi(e))\\
&=h_{r(e)}^{1/2}(\alpha\xi^\prime(\psi(e))+\beta\eta^\prime(\psi(e)))\\
&=\alpha h_{r(e)}^{1/2}\xi^\prime(\psi(e))+\beta h_{r(e)}^{1/2}\eta^\prime(\psi(e))\\
&=(\alpha T_0(\xi^\prime)+\beta T_0(\eta^\prime))(e);
\end{align*}
that is, $T_0$ is a linear map, and hence there exists a lift to a linear map $T:\cE_{Q^\prime}\to\cE_Q$ such that $T\st_{C_c(E^\prime)}=T_0.$
We now claim that $(T,\pi)$ is a coisometric morphism of correspondence $\cE_{Q^\prime}$ to correspondence $\cE_{Q}$.

To see this claim, let $y\in X,$ $e\in E$, $\xi^\prime,\eta^\prime\in C_c(E^\prime)$ and $a^\prime\in C_0(X^\prime)$:
\begin{flalign*}
(1)\; \pi(\innprod{\xi^\prime}{\eta^\prime})(y)&=\innprod{\xi^\prime}{\eta^\prime}(\varphi(y))
=\int\overline{\xi^\prime(e^\prime)}\eta^\prime(e^\prime)\, d\lambda^\prime_{\varphi(y)}(e^\prime)\\
&=\int\overline{(\xi^\prime\circ\psi)(\psi^{-1}(e^\prime))}(\eta^\prime\circ\psi)(\psi^{-1}(e^\prime))h_{r(e)}\,d(\lambda_y\circ\psi^{-1})(e^\prime)\\
&=\int\overline{h_{r(e)}^{1/2}(\xi^\prime\circ\psi)(e_0)}h_{r(e)}^{1/2}(\eta^\prime\circ\psi)(e_0)\, d\lambda_y(e_0)\\
&=\innprod{T(\xi^\prime)}{T(\eta^\prime)}(y)\\\mbox{}\\
(2)\; T(\xi^\prime\cdot a^\prime)(e)&=h_{r(e)}^{1/2}(\xi^\prime\cdot a^\prime)(\psi(e))\\
&=h_{r(e)}^{1/2}(\xi^\prime\circ\psi)(e)(a^\prime(r^\prime(\psi(e))))\\
&=T(\xi^\prime)(e)(a^\prime\circ\varphi)(r(e))\\
&=T(\xi^\prime)\pi(a^\prime)(e)\\
\mbox{and similarly,}\\
(3)\; T(\phi(a^\prime)\xi^\prime)(e)&=\pi(a^\prime)T(\xi^\prime)(e).
\end{flalign*}
Finally, (4) Suppose $a^\prime\in J(\cE_{Q^\prime})=(\phi^\prime)^{-1}(\cK(\cE_{Q^\prime})).$ Then for every $\epsilon>0,$ 
there exist $\theta_\epsilon$ of the form
$$\theta_\epsilon=\sum_{k=1}^{n_\epsilon}\theta_{\xi^\prime_k,\eta^\prime_k}\in \cK(\cE_{Q^\prime})$$ 
where $\xi^\prime_k,\eta^\prime_k\in C_c(E^\prime)$ such that 
$\norm{\theta_\epsilon-\phi^\prime(a^\prime)}<\epsilon.$
To show that $(T,\pi)$ is coisometric, it is enough to show that
$\Psi_T(\phi^\prime(a^\prime))=\phi(\pi(a^\prime))$
for all $a^\prime=\theta_{\xi^\prime,\eta^\prime}$ where $\xi^\prime,\eta^\prime\in C_c(E^\prime).$ To this end, let $\zeta\in C_c(E)$ and $e\in E$. Then
\begin{align*}
\phi(\pi(\theta_{\xi^\prime,\eta^\prime}))\zeta(e)
&=(\theta_{\xi^\prime,\eta^\prime}\circ\varphi)(s(e))\zeta(e)=\theta_{\xi^\prime,\eta^\prime}(s(\psi(e)))\zeta(e)\\
&=(\theta_{\xi^\prime,\eta^\prime}\cdot (\zeta\circ\psi^{-1}))(\psi(e))=\xi^\prime\innprod{\eta^\prime}{\zeta\circ\psi^{-1}}(\psi(e))\\
&=\xi^\prime(\psi(e))\int\overline{\eta^\prime(e^\prime)}(\zeta\circ\psi^{-1})(e^\prime)\,d\lambda^\prime_{r^\prime(\psi(e))}(e^\prime)\\
&=\xi^\prime(\psi(e))\int\overline{(\eta^\prime\circ\psi)(\psi^{-1}(e^\prime))}(\zeta\circ\psi^{-1})(e^\prime)\,d\lambda^\prime_{\varphi(r(e))}(e^\prime)\\
&=\xi^\prime(\psi(e))\int\overline{(\eta^\prime\circ\psi)(\psi^{-1}(e^\prime))}(\zeta\circ\psi^{-1})(e^\prime)\,d\lambda_{r(e)}(\psi^{-1}(e^\prime))\\
&=\xi^\prime(\psi(e))\int\overline{(\eta^\prime\circ\psi)(e_0)}
\zeta(e_0)h_{r(e)}\,d\lambda_{r(e)}(e_0))\\
&=h_{r(e)}^{1/2}(\xi^\prime\circ\psi)(e)\innprod{h_{r(e)}^{1/2}(\eta^\prime\circ\psi)}{\zeta}(r(e))\\
&=T(\xi^\prime)\cdot\innprod{T(\eta^\prime)}{\zeta}(e)\\
&=\phi^\prime(\theta_{T(\xi^\prime),T(\eta^\prime)})\zeta(e)\\
&=\Psi_T(\phi^\prime(\theta_{\xi^\prime,\eta^\prime}))\zeta(e).
\end{align*}
Hence, $(T,\pi)$ is coisometric.

Now take the universal representations $(T_K,\pi_K)$ and $(T_{K^\prime},\pi_{K^\prime})$ that are coisometric on $K=\varphi^\#(K^\prime)=\pi(K^\prime)$ and $K^\prime$, respectively. First, we need to show that $K$ is a closed ideal of $J(\cE_Q).$ Note it is easy to see that $\pi(a)\in\phi^{-1}(\cK(\cE_Q))$ for each $a\in J(\cE_{Q^\prime})$ from the long calculation above and so, $\pi(J(\cE_{Q^\prime}))\subseteq J(\cE_Q)$. 
Furthermore, since $C_c(E)$ is in the image of $T$, we must have $J(\cE_Q)=\pi(J(\cE_{Q^\prime}))$.
If $a\in J(\cE_Q)$ and $b\in K$, then there exist $a^\prime\in J(\cE_{Q^\prime})$ and $b^\prime\in K^\prime$ such that
$a=\pi(a^\prime)$ and $b=\pi(b^\prime)$. So
$$ab=\pi(a^\prime b^\prime)\in \pi(K^\prime)=K.$$
Likewise for $ba\in K$ and closure follows from the fact that $\varphi$ is a homeomorphism.

We now claim that $(T_K,\pi_K)\circ(T,\pi)$ is a morphism that is coisometric on $K^\prime$.  To see this, 
let $\phi_K$ be the left action homomorphism and $a^\prime\in K^\prime$, 
then $\pi(a^\prime)=\varphi^\#(a^\prime)\in K$ and
$$\Psi_{T_K\circ T}(\phi^\prime(a^\prime))=\Psi_{T_K}(\phi(\pi(a^\prime)))=\phi_K(\pi_K\circ\pi(a^\prime)).$$ 

Now by the universal property of $(T_{K^\prime},\pi_{K^\prime}),$ there exists a $*$-homomorphism $\rho:\cO(K^\prime,\cE_{Q^\prime})\to \cO(K,\cE_Q)$ satisfying
$$(T_K,\pi_K)\circ(T,\pi)=\rho\circ (T_{K^\prime},\pi_{K^\prime}).$$
We may do precisely the same with $(T^{-1},\pi^{-1}):\cE_Q\to\cE_{Q^\prime}$ defined by 
$T_0^{-1}(\xi)(e^\prime)=h_{r^\prime(e^\prime)}^{-1/2}(\xi\circ\psi^{-1})(e^\prime)$ for each $e^\prime\in E^\prime$ 
and $\xi\in C_c(E)$ with the lift $T^{-1}:\cE_Q\to\cE_{Q^\prime}$ satisfying $T^{-1}\st_{C_c(E)}=T^{-1}_0.$ Also, $\pi^{-1}(a)=a\circ\varphi^{-1}$ for each $a\in C_0(X)$. Furthermore, note $K^\prime=(\varphi^{-1})^\#(K)$.
Thus, by the universal property of
$(T_K,\pi_K),$ there exists a $*$-homomorphism $\tau:\cO(K,\cE_Q)\to \cO(K^\prime,\cE_{Q^\prime})$ satisfying $$(T_{K^\prime},\pi_{K^\prime})\circ(T^{-1},\pi^{-1})=\tau\circ (T_K,\pi_K)$$
and we note $\rho\circ\tau=id_{\cO(K,\cE_Q)}$ and $\tau\circ\rho=id_{\cO(K^\prime,\cE_{Q^\prime})}.$ 
Hence, $\cO(K,\cE_Q)\cong \cO(K^\prime,\cE_{Q^\prime}).$\\\sq
\end{theorem}

\begin{corollary}\rm With the notation of Theorem \ref{QuiverIso}, if $Q\cong Q^\prime$ then 
$\cT(Q)\cong\cT(Q^\prime)$\index{Toeplitz-Pimsner Algebra} and $C^*(Q)\cong C^*(Q^\prime).$\index{Cuntz-Pimsner Algebra}\\
\pf We need only show $\pi(J_{\cE_{Q^\prime}})=J_{\cE_Q}$. So given any $b\in\ker\phi$, $\xi^\prime\in\cE_{Q^\prime}$ and $e^\prime\in E^\prime,$
\begin{align*}
\phi^\prime(b\circ\varphi^{-1})\xi^\prime(e^\prime)&=b(\varphi^{-1}(s^\prime(e^\prime)))\xi^\prime(e^\prime)\\
&=b(s(\psi^{-1}(e^\prime)))(\xi^\prime\circ\psi)(\psi^{-1}(e^\prime))\\
&=b\cdot (\xi\circ\psi)(\psi^{-1}(e^\prime))=0
\end{align*}
since $\phi(b)=0$, and so $b\circ\varphi^{-1}\in\ker\phi^\prime.$
Furthermore, for each $a^\prime\in J_{\cE_{Q^\prime}},$
$$\phi(\pi(a^\prime))b=(a^\prime\circ\varphi)b=(a(b\circ\varphi^{-1}))\circ\varphi=0.$$
Hence, $\pi(a)\in(\ker\phi^\prime)^\perp$ and recall from the proof of the previous theorem that    $\pi(J(\cE_{Q^\prime}))=J(\cE_Q),$ so $\pi(J_{\cE_{Q^\prime}})\subseteq J_{\cE_Q}$. Conversely, let
$a \in J_{\cE_Q}$. Then there exists $a^\prime\in J(\cE_{Q^\prime})$ such that $\pi(a^\prime)=a$, 
and given any $b^\prime\in\ker\phi^\prime$, 
$$\pi(a^\prime b^\prime)=a\pi(b^\prime)=0$$
since $\pi(b^\prime)\in\ker\phi.$ To see this, let $\xi\in\cE_Q$ and $e\in E.$ Then
\begin{align*} 
\phi(\pi(b^\prime))\xi(e)&=b^\prime(\varphi(s(e)))\xi(e)\\
&=b^\prime(s^\prime(\psi(e)))(\xi\circ\psi^{-1})(\psi(e))\\
&=\phi^\prime(b^\prime)(\xi\circ\psi^{-1})(\psi(e))=0\qquad(\mbox{since $b^\prime\in\ker\phi^\prime$}.)
\end{align*}
Hence, $a^\prime b^\prime=0$ and thus, $a^\prime\in(\ker\phi^\prime)^\perp,$ as desired.\\\sq 
\end{corollary}

\section{Topological Group Quivers}
	\subsection{Definitions}We now investigate a particular class of topological quivers related to a locally compact group and then form and study its associated $C^*$-algebra. This construction leads us to many interesting results including a presentation of the $C^*$-algebra, and the spatial structure of the $C^*$-algebra. 

\begin{definition}\label{TopGrpQuiver}\rm Let $\Gamma$ be a (second countable) locally compact group and let $\alpha,\beta\in\mbox{End}(\Gamma)$ be continuous. Define the closed subgroup, $\Omega_{\alpha,\beta}(\Gamma),$ 
of $\Gamma\times\Gamma,$\index{$\Omega_{\alpha,\beta}(\Gamma)$}
$$\Omega_{\alpha,\beta}(\Gamma)=\{(x,y)\in \Gamma\times\Gamma\st \alpha(y)=\beta(x)\}$$
and let $Q_{\alpha,\beta}(\Gamma)=(\Gamma,\Omega_{\alpha,\beta}(\Gamma), r,s,\lambda)$
\index{Q$\mbox{}_{\alpha,\beta}(\Gamma)$}\index{Topological Quiver}\index{Topological Quiver! $C^*$-algebra Associated with}  where $r$ and $s$ are the group homomorphisms defined by 
$$r(x,y)=x \qquad\mbox{and}\qquad s(x,y)=y$$ 
for each $(x,y)\in\Omega_{\alpha,\beta}(\Gamma)$ and $\lambda_x$ for $x\in\Gamma$ is the measure on 
$$r^{-1}(x)=\{x\}\times\alpha^{-1}(\beta(x))$$ 
defined by 
$$\lambda_x(B)=\mu(y^{-1}s(B\cap r^{-1}(x))\cap\ker\alpha)\qquad\mbox{(for any}\, y\in\alpha^{-1}(\beta(x)))$$
for each measurable $B\subseteq\Omega_{\alpha,\beta}(\Gamma)$ 
where $\mu$ is a left Haar measure (normalized if possible) on $r^{-1}(1_\Gamma)=\{1\}\times\ker\alpha$ 
(a closed normal subgroup of $\Gamma\times\Gamma;$ hence, a locally compact group). Note if $r^{-1}(x)=\emptyset$
then $\alpha^{-1}(\beta(x))=\emptyset$ and so $\lambda_x=0.$ 
This measure is well-defined. If $y_1,y_2\in\alpha^{-1}(\beta(x))$  then
$y_2^{-1}y_1\in\ker\alpha$ and $y_2^{-1}y_1\ker\alpha=\ker\alpha.$ Hence, 
\begin{align*}
\mu(y_1^{-1}s(B\cap r^{-1}(x))\cap\ker\alpha)&=\mu(y_2^{-1}y_1[y_1^{-1}s(B\cap r^{-1}(x))\cap\ker\alpha])\\
&=\mu([y_2^{-1}y_1y_1^{-1}s(B\cap r^{-1}(x))]\cap [y_2^{-1}y_1\ker\alpha])\\
&=\mu(y_2^{-1}s(B\cap r^{-1}(x))\cap\ker\alpha)
\end{align*} for each measurable $B\subseteq \Omega_{\alpha,\beta}(\Gamma).$

Furthermore, let $x\in\Gamma$. If $y\in\alpha^{-1}(\beta(x))\ne\emptyset$ then for any $y_0\in\alpha^{-1}(\beta(x)),$
$$\alpha(y^{-1}y_0)=\alpha(y)^{-1}\alpha(y_0)=\beta(x)^{-1}\beta(x)=1_\Gamma;$$
hence, $y^{-1}y_0\in\ker\alpha$ and $y_0=yy^{-1}y_0\in y\ker\alpha.$ So $\alpha^{-1}(\beta(x))\subseteq y\ker\alpha$
and if $y_0=yz$ for some $z\in\ker\alpha$, then $\alpha(y_0)=\alpha(yz)=\alpha(y)=\beta(x).$ Thus, $\alpha^{-1}(\beta(x))=y\ker\alpha.$ So if $r^{-1}(x)=\emptyset$ then $\mbox{supp }\lambda_x=\emptyset=r^{-1}(x)$ 
and if $(x,y)\in r^{-1}(x)$ then
$$\mbox{supp }\lambda_x=\{x\}\times y\ker\alpha=\{x\}\times \alpha^{-1}(\beta(x))=r^{-1}(x).$$

Moreover, the map $$x\mapsto\int_{r^{-1}(x)}f(e)\,d\lambda_x(e)=\int_{\ker\alpha}f(x,yz)\,d\mu(z)$$
has (closed) support contained in $r(\mbox{supp } f)$ (which is compact for each $f\in C_c(\Omega_{\alpha,\beta}(\Gamma))$.) 
Finally, the combined facts that $\alpha$ and $\beta$ are continuous endomorphisms (hence, open maps) and $f$ has 
compact support guarantees that this map is continuous.

Call $Q_{\alpha,\beta}(\Gamma)$ a \emph{topological group relation.}\index{Topological Group Relation}\index{Topological Relations}
Define $\cE_{\alpha,\beta}(\Gamma)$\index{$\cE_{\alpha,\beta}(\Gamma)$}\index{C$\mbox{}^*$-correspondence} 
to be the $C_0(\Gamma)$-correspondence $\cE_{Q_{\alpha,\beta}(\Gamma)}$ and
form the Cuntz-Pimsner algebra\index{$\cO_{\alpha,\beta}(\Gamma)$! Definition}\index{Cuntz-Pimsner Algebra}
$$\cO_{\alpha,\beta}(\Gamma):=C^*(Q_{\alpha,\beta}(\Gamma))=\cO(J_{\cE_{\alpha,\beta}(\Gamma)},\cE_{\alpha,\beta}(\Gamma))$$
and the Toeplitz-Pimsner algebra\index{$\cT_{\alpha,\beta}(\Gamma)$! Definition}\index{Toeplitz-Pimsner Algebra}
$$\cT_{\alpha,\beta}(\Gamma):=\cT(Q_{\alpha,\beta}(\Gamma)).$$
\end{definition}

\begin{remark}\rm It will be implicitly assumed that $\Gamma$ is second countable. Furthermore, since $\Gamma$ is locally compact Hausdorff, $r^{-1}(x)$ is closed and locally compact. Moreover, whenever $r$ is a local homeomorphism, $r^{-1}(x)$ is discrete and hence, $\lambda_x$ is counting measure (normalized when $\abs{\ker\alpha}<\infty$.)

Further note the nomenclature ``topological group relation'' is chosen because the topological quiver $Q_{\alpha,\beta}(\Gamma)$ is a topological relation as introduced in \cite{BB2}.
\end{remark}

\begin{remark}\rm By Theorem \ref{QuiverIso}, we see that the resulting $C^*$-algebra is independent of the left Haar measure chosen, since any two such left Haar measures will satisfy $\mu=c\mu^\prime$ for some $c>0$. Thus we may solely
regard the normalized left Haar measure when available. 
\end{remark}

\begin{lemma}\rm If $\Gamma$ is a compact group and $\alpha$ is a continuous endomorphism on $\Gamma$, then
$\ker\alpha$ is either finite or an uncountable perfect compact set.\\
\pf Assume $\ker\alpha=\{z_i\}_{i\in I}$ is not finite. Then since $\Gamma$ is compact and $\ker\alpha$ is the preimage of a closed
set, $\ker\alpha$ is compact. Furthermore, let $\lim z_n=z\in\ker\alpha$ where $\{z_n\not=z\}_{n\in\zN}$ is a convergent subsequence. Then, for any $y\in\ker\alpha,$
$$y=\lim_n yz^{-1}z_n$$
but $y\not=yz^{-1}z_n$ for all $n$, so $\ker\alpha$ is a perfect set. The Baire Category Theorem implies any perfect compact 
space must be uncountable.\\\sq
\end{lemma}

We will be mainly concerned with examples where $\ker\alpha$ is finite, but we shall give a non-trivial example when $\abs{\ker\alpha}=\infty$ in a later subsection. If $\abs{\ker\alpha}<\infty$ and $y\in\alpha^{-1}(\beta(x))$ then 
$$\abs{r^{-1}(x)}=\abs{\{x\}\times\alpha^{-1}(\beta(x))}=\abs{\alpha^{-1}(\beta(x))}
=\{y\ker\alpha\}=\abs{\ker\alpha}<\infty.$$ 
Hence, the inner product on $\cE_{\alpha,\beta}(\Gamma)$ becomes
$$\innprod{\xi}{\eta}(x)=\frac{1}{\abs{\ker\alpha}}\sum_{\alpha(y)=\beta(x)}\overline{\xi(x,y)}\eta(x,y).$$
Moreover, one may inquire,``When is $Q_{\alpha,\beta}(\Gamma)$  sourceless/sinkless?'' This is, in fact, quite easy to see.
We note $r(\Omega_{\alpha,\beta}(\Gamma))=\alpha^{-1}(\beta(\Gamma))$ and $s(\Omega_{\alpha,\beta}(\Gamma))=\beta^{-1}(\alpha(\Gamma))$.
Hence, $Q_{\alpha,\beta}(\Gamma)$ is sourceless (or sinkless) if and only if $\alpha^{-1}(\beta(\Gamma))=\Gamma$ (or $\beta^{-1}(\alpha(\Gamma))=\Gamma$.) In particular, $Q_{\alpha,\beta}(\Gamma)$ is sinkless and sourceless if $\alpha$ and $\beta$ are surjective. 

Given a continuous function $\omega:\Gamma\to(0,\infty)$, define a family of Radon measures 
$\{\lambda^\omega_x\st x\in\Gamma\}$ by
$$\lambda^\omega_x(B)=\int_B\omega(s(e))\,d\lambda_x(e)$$
for each measurable $B\subset\Omega_{\alpha,\beta}(\Gamma)$
and let
$$\norm{\xi}_\omega^2:=\norm{\lambda_x^\omega(\xi)}_\infty
=\sup_{x\in\Gamma}\{\int\omega(s(e))\abs{\xi}^2(e)\,d\lambda_x(e)\}$$
for each $\xi\in\cE_{\alpha,\beta}(\Gamma).$

\begin{lemma}\rm Let $\Gamma$ be a compact group. Then $\norm{\cdot}_\omega$ is an equivalent norm of $\norm{\cdot}.$\\
\pf First, $\norm{\cdot}_\omega$ is a seminorm and there exist $m,M>0$ such that 
$0<m\le\omega(x)\le M$ for each $x\in\Gamma$. Thus,
\begin{align*}
\norm{\xi}_\omega^2&=\sup_{x\in\Gamma}\{\int\omega(s(e))\abs{\xi}^2(e)\,d\lambda_x(e)\}\le M\sup_{x\in\Gamma}\{\int\abs{\xi}^2(e)\,d\lambda_x(e)\}=M\norm{\xi}^2
\end{align*} for each $\xi\in\cE_{\alpha,\beta}(\Gamma).$

On the other hand, given $\xi\in\cE_{\alpha,\beta}(\Gamma)$ there exists $x_0\in\Gamma$ 
such that $\norm{\xi}^2=\innprod{\xi}{\xi}(x_0)$. So
\begin{align*}
\norm{\xi}^2_\omega&\ge\lambda_{x_0}^\omega(\xi)=\int\omega(s(e))\abs{\xi}^2(e)\,d\lambda_{x_0}(e)\\
&\ge m\int\abs{\xi}^2(e)\,d\lambda_{x_0}(e)\\
&=m\innprod{\xi}{\xi}(x_0)\\
&=m\norm{\xi}^2.
\end{align*}
Thus, $m\norm{\xi}^2\le\norm{\xi}_\omega^2\le M\norm{\xi}^2$ for each $\xi\in\cE_{\alpha,\beta}(\Gamma).$

If $\norm{\xi}_\omega=0$ then $m\norm{\xi}^2\le0\le M\norm{\xi}^2$. So $\norm{\xi}=0$ and thus, $\xi=0$.
Hence, $\norm{\cdot}_\omega$ is a norm.\\\sq
\end{lemma}

\begin{lemma}\label{Complete}\rm 
Suppose $\Gamma$ is a compact group and $\alpha,\beta\in\mbox{End}(\Gamma)$ are continuous such that 
$\alpha$ has finite kernel. Then the norm given by the inner product, $\norm{\cdot}$, is equivalent to the 
sup-norm and $C(\Omega_{\alpha,\beta}(\Gamma))$ is complete with respect to $\norm{\cdot}.$ Moreover, $\cE_{\alpha,\beta}(\Gamma)=C(\Omega_{\alpha,\beta}(\Gamma))$ is a full correspondence\index{Hilbert $C^*$-module! Full}
\index{C$\mbox{}^*$-correspondence}.\\
\pf Note first that $\Omega_{\alpha,\beta}(\Gamma)$ is a closed subset of the compact set $\Gamma\times\Gamma$, hence $\Omega_{\alpha,\beta}(\Gamma)$ is compact. Note 
$$\abs{r^{-1}(x)}=\begin{cases}
\abs{\ker\alpha}&\mbox{if }r^{-1}(x)\ne\emptyset\\ 0&\mbox{if }r^{-1}(x)=\emptyset
\end{cases}<\infty$$ 
for each $x\in\Gamma$ and so for $\Gamma_0=\{x\in\Gamma\st r^{-1}(x)\ne\emptyset\},$
\begin{align*}
\norm{\xi}^2&=\norm{\innprod{\xi}{\xi}}_\infty\\
&=\max_{x\in\Gamma}\{\int_{r^{-1}(x)}\abs{\xi}(e)\,d\lambda_x(e)\}\\
&=\max_{\Gamma_0}\{\frac{1}{\abs{\ker\alpha}}\sum_{\alpha(y)=\beta(x)} \abs{\xi}^2(x,y)\}\\
&=\frac{1}{\abs{\ker\alpha}}\sum_{z\in\ker\alpha}\max_{x\in\Gamma_0}\{\abs{\xi}^2(x,y_xz)\}\qquad
\mbox{for any $y_x\in\alpha^{-1}(\beta(x))$ (per $x\in\Gamma_0$)}\\
&\le\frac{1}{\abs{\ker\alpha}}\sum_{z\in\ker\alpha}\max_{e\in\Omega_{\alpha,\beta}(\Gamma)}\{\abs{\xi}^2(e)\}\\
&=\norm{\xi}^2_\infty.
\end{align*}
Next, $\frac{1}{\abs{\ker\alpha}}\sum_{z\in\ker\alpha}\abs{\xi}^2(x,yz)\le\norm{\xi}^2$ for each 
$x\in\Gamma$ and $y\in\alpha^{-1}(\beta(x))$. Thus, $\abs{\xi}^2(x,y)\le\abs{\ker\alpha}\norm{\xi}^2$ for each 
$(x,y)\in\Omega_{\alpha,\beta}(\Gamma)$ and so, 
$$\frac{1}{\abs{\ker\alpha}}\norm{\xi}_\infty^2\le\norm{\xi}^2\le\norm{\xi}_\infty^2.$$ 
Since $C(\Omega_{\alpha,\beta}(\Gamma))$ is complete with respect to the sup-norm, we have that $C(\Omega_{\alpha,\beta}(\Gamma))$ is complete with respect to the norm given by the inner product.  

Finally, note that $\cE_{\alpha,\beta}(\Gamma)$ is full since given any $a\in C(\Gamma)$, it is possible to 
define $\xi,\eta\in C(\Omega_{\alpha,\beta}(\Gamma))$
by $\xi(x,y)=1$ and $\eta(x,y)=a(x)$ for $(x,y)\in\Omega_{\alpha,\beta}(\Gamma)$ and note
$$\innprod{\xi}{\eta}(x)=\int \overline{\xi(e)}\eta(e)\,d\lambda_x(e)=a(x).$$\sq
\end{lemma}

\begin{remark}\rm The assertion in the previous lemma concerning the full Hilbert module $\cE_{\alpha,\beta}(\Gamma)$
remains true for infinite kernel. The proof is identical.
\end{remark}

In many cases one can simplify $\alpha$ and $\beta$ and yet preserve the structure of the topological quiver and
hence, the resulting $C^*$-algebra.

\begin{proposition}\label{equiv} \rm Let $\sigma\in\mbox{Aut}(\Gamma)$ be continuous. Then
$$Q_{\alpha,\beta}(\Gamma)\cong Q_{\sigma\circ\alpha, \sigma\circ\beta}(\Gamma)\qquad\mbox{and}\qquad
Q_{\alpha,\beta}(\Gamma)\cong Q_{\alpha\circ\sigma,\beta\circ\sigma}(\Gamma).$$\index{Q$\mbox{}_{\alpha,\beta}(\Gamma)$}
Thus, $$\cO_{\alpha,\beta}(\Gamma)\cong \cO_{\sigma\circ\alpha, \sigma\circ\beta}(\Gamma),\qquad
\cO_{\alpha,\beta}(\Gamma)\cong \cO_{\alpha\circ\sigma,\beta\circ\sigma}(\Gamma),$$\index{$\cO_{\alpha,\beta}(\Gamma)$}
$$\cT_{\alpha,\beta}(\Gamma)\cong \cT_{\sigma\circ\alpha, \sigma\circ\beta}(\Gamma)\qquad\mbox{and}\qquad
\cT_{\alpha,\beta}(\Gamma)\cong \cT_{\alpha\circ\sigma,\beta\circ\sigma}(\Gamma).$$\index{$\cT_{\alpha,\beta}(\Gamma)$}
\pf This is, in fact, a corollary of Theorem \ref{QuiverIso}. For the first isomorphism\index{Topological Quiver! Isomorphism}, use $\varphi=id_\Gamma$ and 
$\psi=id_{\Omega_{\alpha,\beta}(\Gamma)}$ while noting
$$(x,y)\in\Omega_{\alpha,\beta}(\Gamma)\mbox{ if and only if } (x,y)\in\Omega_{\sigma\circ\alpha,\sigma\circ\beta}(\Gamma).$$
Then $\varphi(s(e))=s(\psi(e))$, $\varphi(r(e))=r(\psi(e))$ for all $e\in\Omega_{\alpha,\beta}(\Gamma)$ and
$\lambda_x=\lambda_x^\sigma$ where $\lambda^\sigma$ is the family of Radon measures for $Q_{\sigma\circ\alpha,\sigma\circ\beta}(\Gamma).$

For the second isomorphism, use $\varphi=\sigma^{-1}$ and $\psi(x,y)=(\sigma^{-1}(x),\sigma^{-1}(y))$. Then
$(x,y)\in r^{-1}(x)$ if and only if $\alpha(y)=\beta(x).$ Let $x_0=\sigma^{-1}(x)$ and $y_0=\sigma^{-1}(y)$, then
$(x_0,y_0)\in r_\sigma^{-1}(\varphi(x))$ if and only if $(x,y)\in r^{-1}(x)$ where $r_\sigma$ is the range map of $Q_{\alpha\circ\sigma,\beta\circ\sigma}(\Gamma)$. Furthermore, if $(x,y)\ne (x,y^\prime)$ in $\Omega_{\alpha,\beta}(\Gamma)$,
then $y_0=\sigma^{-1}(y)\ne\sigma^{-1}(y^\prime)=y^\prime_0$ and so, $(x_0,y_0)\ne (x_0,y_0^\prime)$ in $\Omega_{\alpha\circ\sigma,\beta\circ\sigma}(\Gamma)$. Thus, we have shown
$$\lambda_{\varphi(x)}^\sigma\circ\psi=c\lambda_x \mbox{ for some $c>0$}$$
where $\lambda^\sigma$ is the family of Radon measure for $Q_{\alpha\circ\sigma,\beta\circ\sigma}(\Gamma)$ 
since both are left Haar measures on the same set.\\
Also, for $e=(x,y)\in\Omega_{\alpha,\beta}(\Gamma)$
$$\varphi(s(e))=\sigma^{-1}(y)=s_\sigma(\psi(e))\qquad\mbox{and}\qquad\varphi(r(e))=\sigma^{-1}(x)=r_\sigma(\psi(e)).$$
\sq
\end{proposition}

\begin{example}\cite{EaHR} \rm Let $\alpha$ be an endomorphism of a unital $C^*$-algebra $B$. 
A \emph{transfer operator}\index{Transfer Operator} $L$ for $(B, \alpha)$ is a positive linear map $L:B\to B$ such that $L(\alpha(a)b) = aL(b)$ for all $a, b\in B$. 
\cite{EaHR} calls the triple $(B, \alpha, L)$ an Exel system and create the Cuntz-Pimsner algebra \index{Cuntz-Pimsner Algebra}
out of the right Hilbert module $M_L,$ the completion of $B$ with respect to the left and right actions
$$a\cdot\xi\cdot b=a\xi\alpha(b)$$
and inner product
$$\innprod{\xi}{\eta}=L(\xi^*\eta).$$
When $B=C(\Gamma)$ for compact group $\Gamma$ and $\alpha=\beta^\#$ where $\beta\in\mbox{End}(\Gamma)$, we get 
$$L(\xi)(x)=\lambda_x(\xi)$$
defines a transfer operator for $(C(\Gamma),\alpha)$ where $\lambda$ is the family of Radon measures 
for $Q_{\beta,1}(\Gamma).$
First note that $\Omega_{\beta,1}(\Gamma)=\{(x,y)\st x=\beta(y)\}$ and so $C(\Omega_{\beta,1}(\Gamma))\cong C(\Gamma)$ where the left and right actions become
$$a\cdot\xi\cdot b=a\xi\beta^\#(b)$$
and inner product
$$\innprod{\xi}{\eta}(x)=\int_{s(r^{-1}(x))}\overline{\xi(y)}\eta(y)\,d\lambda_x(x,y)=L(\xi^*\eta)(x)$$
for $a,b\in C(\Gamma),$ $\xi,\eta\in M_L=\cE_{\beta,1}(\Gamma)$ and $x\in\Gamma.$ 
Now for $x\in\Gamma,$ $a\in C(\Gamma)$ and $\xi\in M_L$,
\begin{align*}
L(\alpha(a)\xi)(x)&=\lambda_x(\alpha(a)\xi)\\
&=\int a(\beta(y))\xi(y)\,d\lambda_x(x,y)\\
&=\int a(x)\xi(y)\,d\lambda_x(x,y)\\
&=a(x)\lambda_x(\xi)\\
&=[aL(\xi)](x).
\end{align*}
This shows the $C^*$-algebras $\cO_{\beta,1}(\Gamma)$ were also studied in \cite{EaHR}. Moreover,
for the $C^*$-algebras $\cO_{\alpha,\beta}(\Gamma),$ notice
\begin{align*}
\lambda_x(\alpha^\#(a)\cdot\xi)&=\int a(\alpha(s(e)))\xi(e)\,d\lambda_x(e)\\
&=\int a(\beta(r(e)))\xi(e)\,d\lambda_x(e)\\
&=\int a(\beta(x))\xi(e)\,d\lambda_x(e)\\
&=\beta^\#(a)(x)\lambda_x(\xi)
\end{align*}
and hence, these topological group relations\index{Topological Group Relation} extend the context in \cite{EaHR}

However, there are topological group relations that are not of this form; for instance, $\cO_{1,0}(\zZ_3)\not\cong\cO_{n,1}(\zZ_3)$ for any $n=0,1,2$ (see  Example 3.12.)
\end{example}

\begin{example}\rm Let $\alpha\in\mbox{End }(\Gamma)$ be surjective for some compact group $\Gamma$. Then 
$$\Omega_{1,\alpha}(\Gamma)=\{(x,\alpha(x))\st x\in \Gamma\}\cong \Gamma.$$\index{$\Omega_{\alpha,\beta}(\Gamma)$}
Hence, the left and right actions become
$$(a\cdot \xi\cdot b)(x,\alpha(x))=a(\alpha(x))\xi(x,\alpha(x))b(x)$$
and the inner product becomes
$$\innprod{\xi}{\eta}(x)=\overline{\xi(x,\alpha(x))}\eta(x,\alpha(x))$$
for $\xi,\eta\in C(\Omega_{1,\alpha}(\Gamma)),$ $a,b\in C(\Gamma)$ and $x\in\Gamma.$
We quickly note that, for $\phi$ the left action of $C(\Gamma)$, $\ker\phi=\{0\}$ since $\alpha$ is surjective and 
so $J_{\cE_{1,\alpha}(\Gamma)}=J(\cE_{1,\alpha}(\Gamma))=\phi^{-1}(\cK(C(\Omega_{1,\alpha}(\Gamma)))).$
\index{J$\mbox{}_\cE$}\index{J$(\cE)$}\index{$\cE_{\alpha,\beta}(\Gamma)$}\index{Q$\mbox{}_{\alpha,\beta}(\Gamma)$}
\index{$\cO_{\alpha,\beta}(\Gamma)$}
We also note, for any $a\in C(\Gamma)$, $\xi(x,y)=a(y)$ defines $\xi\in C(\Omega_{1,\alpha}(\Gamma))$ and so
$$\theta_{\xi,1}(\eta)(x,\alpha(x))=\xi(x,\alpha(x))\innprod{1}{\eta}(x)=a(\alpha(x))\eta(x,\alpha(x))=a\cdot\eta (x,\alpha(x))$$
for all $(x,\alpha(x))\in\Omega_{1,\alpha}(\Gamma)$ and $\eta\in C(\Omega_{1,\alpha}(\Gamma));$ that is,
$$J_{\cE_{1,\alpha}(\Gamma)}=C(\Gamma).$$

Now we let $u\in C(\Gamma)$ be defined by $u(x,\alpha(x))=1$ for all $(x,\alpha(x))\in\Omega_{1,\alpha}(\Gamma).$
Then note
$$[a\cdot u](x,\alpha(x))=a(\alpha(x))=[u\cdot \alpha^\#(a)] (x,\alpha(x)),$$
$$\innprod{u}{u}=1\qquad\mbox{and}\qquad u\innprod{u}{\xi}=\xi$$
where $\alpha^\#(a)(x)=a(\alpha(x))$.
Hence, by Theorem \ref{BasisGen} with $S=S_u$,
$$\cO_{1,\alpha}(\Gamma)\cong C^*(C(\Gamma), S\st S^*S=SS^*=1\mbox{ and }\alpha^\#(a)=S^*aS\mbox{ for all $a\in C(\Gamma)$}).$$
That is, if $\alpha^\#$ is an automorphism, then
$$\cO_{1,\alpha}(\Gamma)\cong C(\Gamma)\rtimes_{\alpha^\#}\zZ$$
as (indirectly) noted by Pimsner \cite{Pims}
and in particular, if $\alpha\in\mbox{Aut}(\Gamma)$ then
$$\cO_{1,\alpha}(\Gamma)\cong C(\Gamma)\rtimes_{\alpha^\#}\zZ\qquad\mbox{and}\qquad\cO_{\alpha,1}(\Gamma)\cong\cO_{1,\alpha^{-1}}(\Gamma)\cong C(\Gamma)\rtimes_{(\alpha^{-1})^\#}\zZ.$$
\end{example}

	\subsection{Topological $\zZ_p-$Quivers}\begin{example}\rm Let $\Gamma$ be the compact abelian group $\zZ_p=\zZ/p\zZ$.\index{Topological Group Relation}
\index{Topological Group Relation! $Z_p$-quiver}\index{Topological Relations}\index{Graph $C^*$-algebra}
Then $\cO_{\alpha,\beta}(\Gamma)$\index{$\cO_{\alpha,\beta}(\Gamma)$} 
is, in general, a graph $C^*$-algebra, but this does not yield all such since 
$(1_\Gamma,1_\Gamma)\in \Omega_{\alpha,\beta}(\Gamma)$ for all $\alpha$ and $\beta$. Recall $\mbox{End}(\zZ_p)\cong \zZ_p$.
\end{example}

\begin{figure}[htb]
\center{\includegraphics[height=1cm]{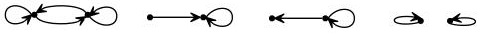}}
\caption[The $\zZ_2$-quivers]{\label{fig:4} The $\zZ_2$-quivers: From left to right: $Q_{0,0}(\zZ_2)$, $Q_{0,1}(\zZ_2)$, 
$Q_{1,0}(\zZ_2)$, $Q_{1,1}(\zZ_2)$}
\end{figure}
\begin{figure}[htb]
\center{\includegraphics[height=3cm]{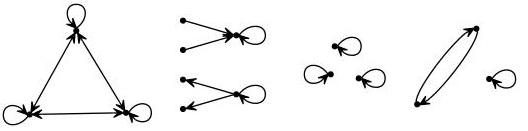}}
\caption[The $\zZ_3$-quivers]{\label{fig:5} The $\zZ_3$-quivers: 
 From left to right and top to bottom: $Q_{0,0}(\zZ_3)$, $Q_{0,1}(\zZ_3)\cong Q_{0,2}(\zZ_3)$, 
$Q_{1,0}(\zZ_3)\cong Q_{2,0}(\zZ_3)$, $Q_{1,1}(\zZ_3)\cong Q_{2,2}(\zZ_3)$, $Q_{1,2}(\zZ_3)\cong Q_{2,1}(\zZ_3)$}
\end{figure}

\begin{example}\rm It can be shown that the graph $C^*$-algebra of a graph with a single loop on one vertex is $\cO_1\cong C(\zT).$ Hence, $\cO_{1,1}(\zZ_2)\cong C(\zT)\oplus C(\zT),$ and 
$\cO_{1,1}(\zZ_3))\cong C(\zT)^{\oplus 3}.$ Furthermore, it is shown in \cite{HMS} that 
$\cO_{0,1}(\zZ_3)\cong C(S^2_{q,1})$ for any $q\in\zR\setminus[0,1]$ (see \cite{HMS} for further details).
Finally, the author shows $\cO_{1,0}(\zZ_2)\cong\cT$ and $\cO_{1,2}(\zZ_3)=C(\zT)\oplus M_2(C(\zT))$ in \cite{McThesis}.
\end{example}

Just a few interesting facts concerning these types of topological group quivers follow: 

\begin{proposition}\rm Let $n\in\zZ_p$ be non-zero. Then
$$Q_{n,n}(\zZ_p)\cong\sqcup_{j=1}^{p/\gcd(n,p)}Q_{0,0}(\zZ_{\gcd(n,p)}).$$
Moreover, 
$$\cO_{n,n}(\zZ_p)\cong\oplus_{j=1}^{p/\gcd(n,p)}\cO_{0,0}(\zZ_{\gcd(n,p)}).$$
\pf Let $z=e^{2\pi i/p}$ and let $n_0=\gcd(n,p).$ Denote $Z_k=\{e^{2(k+mp/n_0)\pi i/p}\st m=0,...,n_0-1\}$ for
$k=0,...,p/n_0-1.$ Then note for any $x,y\in Z_k$ where $x=e^{2(k+m_1p/n_0)\pi i/p}$ and $y=e^{2(k+m_2p/n_0)\pi i/p},$
$$ (xy^{-1})^{n_0}=e^{2(m_1-m_2)\pi i}=1;$$
that is, $x^{n_0}=y^{n_0}.$
Furthermore, if $x=e^{2(k_1+m_1p/n_0)\pi i/p}\in Z_{k_1}$ and $y=e^{2(k_2+m_2p/n_0)\pi i/p}\in Z_{k_2}$  for $k_1\ne k_2$ 
then
$$(xy^{-1})^{n_0}=e^{2(k_1-k_2)n_0\pi i/p + 2(m_1-m_2)\pi i}=e^{2(k_1-k_2)n_0\pi i/p}$$
but $\abs{k_1-k_2}< p/n_0$ ensuring that $(xy^{-1})^{n_0}\ne 1$ (or better written, $x^{n_0}\ne y^{n_0}$.)
The rest is now a trivial application of Theorem \ref{QuiverIso}.\\\sq
\end{proposition}

It can be shown (see \cite{Kum}, \cite{Ev}, \cite{aHR}) that the graph 
with $m$-points connected in a single length $m$ loop admits the universal $C^*$-algebra
$M_m(C(\zT)).$ Let $\zZ_p^*$ denote the multiplicative group $\{n\in\zZ_p\st \gcd(n,p)=1\}.$ 

\begin{proposition}\rm Let $n,p\in\zN$ such that $\gcd(n,p)=1$. Then:
\begin{enumerate}
\item $Q_{n,1}(\zZ_p)$ consists of disjoint loops of various lengths less than or equal to 
$o_p(n)$ (the order of $n$ in $\zZ_p^*$).
\item $Q_{1,n}(\zZ_p)\cong Q_{n,1}(\zZ_p).$
\item The number of base points of loops of length $k\le o_p(n)$
\index{Loop}\index{Base Point} in $Q_{n,1}(\zZ_p)$ is $$\gcd(n^k-1,p).$$ 
\end{enumerate}
\pf (1) Given $z\in\zZ_p$, there exists a unique edge with range $z$, likewise source $z$, since $\gcd(n,p)=\gcd(1,p)=1.$
Furthermore, there exists a number $k\le o_p(n)$ such that $z^{n^k}=z.$

(2) It is clear from (1) that reversing the arrows will not change the graph ($Q_{1,n}(\zZ_p)$ is $Q_{n,1}(\zZ_p)$ with reversed arrows).

(3) We need only find all solutions, $y$, to $n^ky=y \mod p.$ This equation has $\gcd(n^k-1,p)$ distinct solutions.\\\sq 
\end{proposition}

\begin{example}\rm With this proposition, it is easy to determine the $C^*$-algebra $\cO_{n,1}(\zZ_p)$ where $\gcd(n,p)=1.$
\begin{enumerate}
\item $\cO_{5,1}(\zZ_{72}):$ Note $o_{72}(5)=6$ and
\begin{align*}
\gcd(5-1,72)&=\gcd(4,72)=4\\
\gcd(5^2-1,72)&=\gcd(24,72)=24\\
\gcd(5^3-1,72)&=\gcd(124,72)=4\\
\gcd(5^4-1,72)&=\gcd(624,72)=24\\
\gcd(5^5-1,72)&=\gcd(3124,72)=4\\
\gcd(5^6-1,72)&=\gcd(15624,72)=72.
\end{align*}
Hence, there are 4 base points of 1-loops and 24 base points of 2-loops. Note: each 1-loop counts here, and so there are
20 base points of 2-loops that aren't base points of 1-loops; that is, 10 loops of length 2.
Furthermore, the 4 base points of 3-loops are already counted by the 1-loops, and so on... 
Finally, there are 72 base points of 6-loops. Subtracting all possible base points of 2-loops, we obtain 48 base points of 6-loops
and hence, 8 6-loops. Therefore,
$$\cO_{5,1}(\zZ_{72})=C(\zT)^{\oplus 4}\oplus M_2(C(\zT))^{\oplus 10}\oplus M_6(C(\zT))^{\oplus 8}.$$
Remark: $4+2(10)+6(8)=72.$     
\item $\cO_{6,1}(\zZ_{77}):$ Note $o_{77}(6)=10$ and
$$\gcd(6-1,77)=1\qquad \gcd(6^2-1,77)=7\qquad\gcd(6^{10}-1,77)=77.$$
That is, there is one copy of $C(\zT)$, $\frac{7-1}{2}=3$ copies of $M_2(C(\zT))$, and $\frac{77-7}{10}=7$ copies of 
$M_{10}(C(\zT));$ i.e.,
$$\cO_{6,1}(\zZ_{77})=C(\zT)\oplus M_2(C(\zT))^{\oplus 3}\oplus M_{10}(C(\zT))^{\oplus 7}.$$
Remark: $1+2(3)+10(7)=77.$
\end{enumerate}
\end{example}

\begin{proposition}\rm Let $p$ be prime and let $\varphi(p-1)$ denote the number of divisors of $p-1$. Then the number of 
distinct (non-isomorphic) topological $\zZ_p$-quivers is $\varphi(p-1)+3$.\\
\pf Begin by noting that $\mbox{Aut}(\zZ_p)=\zZ_p^*=\{1,...,p-1\}$, the numbers coprime to the prime number $p$.
Hence, by Proposition \ref{equiv}, 
$$Q_{n,m}(\zZ_p)\cong Q_{kn,km}(\zZ_p)$$
for each $k=1,...,p-1$. This gives us only the following quivers to check
$$Q_{0,0}(\zZ_p), Q_{0,1}(\zZ_p), Q_{1,0}(\zZ_p)\mbox{ [the $3$ in $\varphi(p-1)+3$], and } Q_{1,n}(\zZ_p)$$
where $n=1,...,p-1$. We see immediately that $Q_{0,0}(\zZ_p),Q_{0,1}(\zZ_p)$ and $Q_{1,0}(\zZ_p)$ are all distinct since
$\Omega_{0,0}(\zZ_p)=\zZ_p\times\zZ_p$, $\Omega_{0,1}(\zZ_p)=\{1\}\times\zZ_p$ (i.e., all arrows have range equal 1) and
$\Omega_{1,0}(\zZ_p)=\zZ_p\times\{1\}$ (i.e., all arrows have source equal 1.) Furthermore, $Q_{1,n}(\zZ_p)$ has
one and only one arrow with range $k$ (likewise, with source $k$) so $Q_{1,n}(\zZ_p)$ is not isomorphic to any of $Q_{0,0}(\zZ_p), Q_{0,1}(\zZ_p),$ or $Q_{1,0}(\zZ_p).$ In fact, $Q_{1,n}(\zZ_p)$ is solely defined based on the length of its loops; that is,
the order of $n$ in $\zZ_p^*$ as the previous proposition proves. Since all orders divide $p-1$ and all orders are achieved, 
the resulting count of distinct topological $\zZ_p$-quivers is $\varphi(p-1)+3.$\\\sq 
\end{proposition}

	\subsection{Topological $\zT^d$-Quivers}\begin{example}\label{ToriQuiver}\rm\index{Topological Group Relation! $\zT^d$-quivers}\index{Topological Group Relation}
\index{Topological Relations}\index{Topological Quiver}\index{$\cO_{F,G}(\zT^d)$! Definition}
For the compact abelian group $\zT^d,$ note $\mbox{End}(\zT^d)\cong M_d(\zZ)$ (\cite{W}); that is, 
an element $\sigma\in\mbox{End}(\zT^d)$ is of the form $\sigma_F$ for some $F\in M_d(\zZ)$ where 
$$\sigma_F(e^{2\pi it})=e^{2\pi i Ft}$$\index{$\sigma_F$}
for each $t\in\zZ^d.$ To simplify notation, use $F$ and $G$ in place of $\sigma_F$ and $\sigma_G$ whenever convenient.
For instance,
$$Q_{F,G}(\zT^d):=Q_{\sigma_{F},\sigma_{G}}(\zT^d)$$\index{Q$\mbox{}_{\alpha,\beta}(\Gamma)$} 
and the $C^*$-correspondence
$$\cE_{F,G}(\zT^d):=\cE_{\sigma_F,\sigma_G}(\zT^d)$$
where $F,G\in M_d(\zZ)$. We will consider the cases when these maps are surjective; that is,
$\det F$ and $\det G$ are non-zero.

Let $F, G\in M_d(\zZ)$ where $\det F,\det G\ne 0$. Then $\abs{\ker\sigma_F}=\abs{\det F}$
and so, the $C(\zT^d)$-valued inner product becomes
$$\innprod{\xi}{\eta}(x)=\frac{1}{\abs{\det F}}\sum_{\sigma_{F}(y)=\sigma_{G}(x)}\overline{\xi(x,y)}\eta(x,y)$$
for $\xi,\eta\in\cE_{F,G}(\zT^d)=C(\Omega_{F,G}(\zT^d))$ [by Lemma \ref{Complete}] and $x\in\zT^d.$
This is a finite sum since the number of solutions, $y,$ 
to $\sigma_F(y)=\sigma_G(x)$ given any $x\in\zT^d$ is $\abs{\det F}<\infty.$
\end{example}

\begin{remark}\rm The left action, $\phi,$ is defined by
$$\phi(a)\xi(x,y)=a(y)\xi(x,y)$$
for $a\in C(\zT^d),$ $\xi\in C(\Omega_{F,G}(\zT^d))$ and $(x,y)\in\Omega_{F,G}(\zT^d).$  
Note: $\phi$ is injective. To see this claim, let $a\in C(\zT^d)$ and assume $\phi(a)\xi=0$ for each $\xi\in C(\Omega_{F,G}(\zT^d)).$ Then $a(y)\xi(x,y)=0$ for each $(x,y)\in\Omega_{F,G}(\zT^d)$ and $\xi\in C(\Omega_{F,G}(\zT^d)).$ Since 
$s(\Omega_{F,G}(\zT^d))=\{y\in\zT^d\st (x,y)\in\Omega_{F,G}(\zT^d)\}=\zT^d$ by the surjectivity of $\sigma_F,$
$a=0.$ 
\end{remark}

Examples of the form $G=1_d$ are studied in \cite{EaHR} and examples when $d=1$ can be found in \cite{Yam} where $\gcd(n,m)=1.$ Brenken \cite{BB2} examines cases when $F=a1_d$ and $G$ is arbitary, and examines the associated $C^*$-algebra when $d=1$, $F=a\in\zN$ and $G=1$.

In fact, a basis for $\cE_{F,G}(\zT^d)$ can be provided using the following lemma and corollaries to Proposition \ref{equiv}.

\begin{proposition}\label{equiv1}\rm If $U,F,G\in M_d(\zZ)$ such that $\det U=1$ and $\det F,\det G\ne0$, then
$\cO_{F,G}(\zT^d)\cong\cO_{UF,UG}(\zT^d)\cong\cO_{FU,GU}(\zT^d).$\\
\pf This is evident from Proposition \ref{equiv}.\sq
\end{proposition}

\begin{corollary}\label{Smith}\rm 
Given $F,G\in M_d(\zZ)$ such that $\det F, \det G\ne0$, there exist $F^\prime,G^\prime\in M_d(\zZ)$
with $F^\prime$ a positive diagonal matrix and 
$$\cO_{F,G}(\zT^d)\cong\cO_{F^\prime,G^\prime}(\zT^d).$$
\pf Using the Smith normal form (see \cite{DF}), there exist unimodular matrices $U, V\in M_d(\zZ)$ 
and a positive diagonal matrix $D\in M_d(\zZ)$
such that $$F=UDV.$$
Hence $U^{-1}, V^{-1}\in M_d(\zZ)$ are unimodular matrices and, by the previous proposition,
$$\cO_{F,G}(\zT^d)\cong\cO_{D,U^{-1}GV^{-1}}(\zT^d).$$\sq 
\end{corollary}

The following lemma is well known:

\begin{lemma}\rm Given $\Gamma_0$ a finite subgroup of $\zT^d$, $\sum_{z\in\Gamma_0} z=0$ (in $\zC^d$.)
\end{lemma}
 
Let $F=\diag{a_1,...,a_d}\in M_d(\zZ), G=(b_{jk})_{j,k=1}^d\in M_d(\zZ)$ where $a_j>0$ for each $j=1,...,d,$ $\det G\ne0$ and
let $G_j$ denote the $j$-th row of $G$, $(b_{jk})_{k=1}^d.$ Further, let $N=\det F=\prod_{j=1}^d a_j>0$ and let $$\fI(F)=\{\nu=(\nu_j)_{j=1}^d\in\zZ^d\st 0\le \nu_j\le a_j-1\}.\index{$\fI(F)$}$$ 
The $C(\zT^d)$-valued  inner product becomes
$$\innprod{\xi}{\eta}(x)=\frac{1}{N}\sum_{\sigma_F(y)=\sigma_G(x)} \overline{\xi(x,y)}\eta(x,y)$$
for all $\xi,\eta\in C(\Omega_{F,G}(\zT^d))$ and $x\in \zT^d.$

Given $\nu\in\fI(F)$, define $u_\nu\in C(\Omega_{F,G}(\zT^d))$ by 
$$u_\nu(x,y)=y^\nu=\prod_{j=1}^d y^{\nu_j}$$\index{Basis! Orthonormal}
for $(x,y)\in \Omega_{F,G}(\zT^d).$ So $u_{(0,0,...,0)}=1\in C(\Omega_{F,G}(\zT^d)).$
Then
\begin{align*}
\innprod{u_\nu}{u_{\nu^\prime}}(x)&=\frac{1}{N}\sum_{\sigma_F(y)=\sigma_G(x)}\overline y^\nu y^{\nu^\prime}\\
&=\frac{1}{N}\sum_{z\in\ker\sigma_F}z_0^{\nu-\nu^\prime}z^{\nu-\nu^\prime}
\qquad\mbox{for any $z_0\in\sigma_F^{-1}(\sigma_G(x))$}\\
&=\begin{cases}0 &\mbox{if }\nu\ne\nu^\prime\\ 1&\mbox{if }\nu=\nu^\prime\end{cases}=\delta_{\nu}^{\nu^\prime}.
\end{align*}
Furthermore, since
$$u_\nu\innprod{u_\nu}{\xi}(x,y)=u_\nu(x,y)\innprod{u_\nu}{\xi}(x)=y^\nu\innprod{u_\nu}{\xi}(x),$$
we have
\begin{align*}
[\sum_{\nu\in\fI(F)} u_\nu\innprod{u_\nu}{\xi}](x,y)
&=\frac{1}{N}\sum_{\nu\in\fI(F)} \sum_{\sigma_F(w)=\sigma_G(x)}y^\nu\overline{w}^\nu\xi(x,w)\\
&=\frac{1}{N}\sum_{\nu}\sum_{\sigma_F(w)=\sigma_G(x)}[\prod_{j=1}^d y_j^{\nu_j}\overline w_j^{\nu_j}]\xi(x,w).\\
\end{align*}
Note: if $x=e^{2\pi it}$ and $\sigma_F(w)=\sigma_G(x),$ then $w_j^{a_j}=e^{2\pi i(Gt)_j}$ where $(Gt)_j$ denotes the 
$j$-th coordinate of the vector $Gt.$ 
Now note
$$\sum_{k=0}^{a_j-1} y_j^k\overline w_j^k= 
\begin{cases}
\frac{(y_j\overline w_j)^{a_j}-1}{y_j\overline w_j -1} & \text{if $y_j\ne w_j$} \\
a_j & \text{if $y_j=w_j$} \\
\end{cases}
=\begin{cases}
0 & \text{if $x_j\ne w_j$} \\
a_j & \text{if $x_j=w_j$} \\
\end{cases}$$
since $y_j^{a_j}\overline w_j^{a_j}=e^{2\pi i (Gt)_j}e^{-2\pi i (Gt)_j}=1$ where $x=e^{2\pi it}$ $(t\in\zZ^d).$ Hence,
\begin{align*}
[\sum_{\nu\in\fI(F)} u_\nu\innprod{u_\nu}{\xi}](x,y)
&=\frac{1}{N}\sum_{\nu}\sum_{\sigma_F(w)=\sigma_G(x)}[\prod_{j=1}^d y_j^{\nu_j}\overline w_j^{\nu_j}]\xi(x,w)\\
&=\frac{1}{N}[\prod_{j=1}^d a_j]\xi(x,y)=\xi(x,y).
\end{align*}
Thus, $\{u_\nu\}_{\nu\in\fI(F)}$ is a basis for $\cE_{F,G}(\zT^d)$.

For $\nu\in\fI(F)$ and $j\in\{1,...,d\}$, using notation from Theorem \ref{BasisGen}, set $S_\nu=S_{u_\nu},$ $S=S_{(0,0,...,0)}$ and let $\{U_j\}_{j=1}^d$ be the full spectrum unitaries in $C(\zT^d)$ defined by 
$U_j(y)=y_j$ for $y\in\zT^d.$ Then, using the notation 
$$U^\nu=U_1^{\nu_1}U_2^{\nu_2}\cdot\cdot\cdot U_d^{\nu_d},$$
we have
$$U^\nu S = S_{U^\nu\cdot u_0}=S_\nu$$
since $[U^\nu\cdot u_0](x,y)=y^\nu=u_\nu(x,y),$ and
$$U_j^{a_j} S=S_{U_j^{a_j}\cdot u_0}=S\prod_{k=1}^dU_k^{b_{jk}}=SU^{G_j}$$
since 
$[U_j^{a_j}\cdot u_0](x,y)=y_j^{a_j}=u_0\cdot U_j\circ\sigma_G(x)=u_0\cdot\prod_{k=1}^dU_k^{b_{jk}}(x,y)=u_0\cdot U^{G_j}(x,y).$ 

\begin{remark}\label{FullSpec} \rm \cite[Proposition 2.21]{MS0} states that if the left action $\phi$ is injective then the natural inclusion $C(\zT^d)\hookrightarrow\cO_{F,G}(\zT^d)$ is injective. Hence, the image of $U_j$ in $\cO_{F,G}(\zT^d)$ (still denoted $U_j$) is a full spectrum unitary since $U_j$ is a full spectrum unitary in $C(\zT^d).$
\end{remark} 

Thus, by Theorem \ref{BasisGen}:
\begin{theorem}\label{TdQuiverGen}\rm \index{$\cO_{F,G}(\zT^d)$}\index{$\cO_{\alpha,\beta}(\Gamma)$! Presentation}
\index{Cuntz-Pimsner Algebra}
Let $F=\diag{a_1,...,a_d}, G\in M_d(\zZ)$ where $\det F, \det G\ne0$ and let $G_j$ be the $j$-th row vector of $G$. Further, let $\fI(F)$ denote the set $\{\nu=(\nu_j)_{j=1}^d\in\zZ^d\st 0\le \nu_j\le a_j-1\}$. Then
$\cO_{F,G}(\zT^d)$ is the universal $C^*$-algebra generated by isometries $\{S_\nu\}_{\nu\in\fI(F)}$ and (full spectrum) commuting unitaries $\{U_j\}_{j=1}^d$ that satisfy the relations
\begin{enumerate}
\item $S_\nu^*S_{\nu^\prime}=\innprod{u_\nu}{u_{\nu^\prime}}=\delta_{\nu}^{\nu^\prime},$
\item $U^\nu S=S_\nu$ for all $\nu\in\fI(F),$
\item $U_j^{a_j}S=SU^{G_j},$ for all $j=1,...,d$ and
\item $1=\sum_{\nu\in\fI(F)} S_\nu S_\nu^*=\sum_{\nu\in\fI(F)} U^\nu SS^*U^{-\nu}$
\end{enumerate}
where $U^\nu$ denotes $\prod_{j=1}^dU_j^{\nu_j}.$
\end{theorem}

\begin{corollary}\rm Let $n\in\zN$ and $G\in M_d(\zZ)$ where $\det G\ne 0$, then
$\cO_{n,G}(\zT^d):=\cO_{n1_d, G}(\zT^d)$ is the universal $C^*$-algebra
$$C^*(\{U_j\}_{j=1}^d, S\st U_j^*U_j=U_jU_j^*=\sum_{\nu\in\mathfrak n^d} U^\nu SS^*U^{-\nu}=1, S^*U^\mu S=\delta_0^\mu\;(\mu\in\mathfrak n^d), U_j^nS=SU^{G_j})$$
where $\mathfrak n^d=\{0,...,n-1\}^d.$
\end{corollary}

\begin{corollary}\rm Let $n\in\zN$ and $m\in\zZ\setminus\{0\},$ then
$\cO_{n,m}(\zT)$ is the universal $C^*$-algebra
$$C^*(U, S\st U^*U=UU^*=\sum_{k=0}^{n-1} U^k SS^*U^{-k}=1, S^*U^k S=\delta_0^k\;(0\le k\le n-1), U^nS=SU^m).$$
\end{corollary}
	\subsection{A Topological $\zT^\omega$-Quiver and Some Generalizations}We describe an example with a perfect compact $\ker\alpha$ which is $\zT.$ 
\index{Topological Group Relation! $\zT^\omega$-quiver} 
Let $\Gamma=\zT^\omega:=\prod_{n\in\zN} \zT$ be the compact abelian topological group with the product topology. Tychonoff's Theorem guarantees that $\Gamma$ is a compact group. Let $\sigma\in\mbox{End}(\zT^\omega)$ be the (surjective) group endomorphism
$$\sigma(t_1,t_2,t_3,...)=(t_2,t_3,...)$$
for each $(t_1,t_2,t_3,...)\in\zT^\omega.$
Then 
$$\Omega_{\sigma, 1}(\zT^\omega)=\{(x,y)\in\zT^\omega\times\zT^\omega\st \sigma(y)=x\},$$
\index{$\Omega_{\alpha,\beta}(\Gamma)$}
$$r^{-1}(x)=\{x\}\times\{(z,x)\st z\in\zT\},$$
where $(z,x)=(z,x_1,x_2,...)\in\zT^\omega,$ and finally, $\lambda_x$ is the Haar measure on $\zT.$ 
The $C(\zT^\omega)$-valued inner product becomes
$$\innprod{\xi}{\eta}(x)=\int_\zT\overline{\xi(x,(z,x))}\eta(x,(z,x))\, dz$$
for each $x\in\zT^\omega.$

Also note, given $a\in C(\zT^\omega)$ and $\xi\in C(\Omega_{\sigma,1}(\zT^\omega))$,
$$\sigma^\#(a)\cdot\xi(x,y)=a(\sigma(y))\xi(x,y)=a(x)\xi(x,y)=(\xi\cdot a)(x,y)$$
for each $(x,y)\in\Omega_{\sigma,1}(\zT^\omega).$ Hence,
$$\sigma^\#(a)S_\xi=S_{\sigma^\#(a)\cdot\xi}=S_{\xi\cdot a}=S_\xi a$$
for each $a\in C(\zT^\omega)$ and $\xi\in C(\Omega_{\sigma,1}(\zT^\omega))$ where $S_\xi$ is the image of $\xi$
into the universal $C^*$-algebra $\cO_{\sigma,1}(\zT^\omega).$

This result is, in fact, a special case of the following:

\begin{proposition}\rm Let $\Gamma$ be a locally compact group and $\alpha,\beta\in\mbox{End}(\Gamma).$ For any 
$\xi\in C(\Omega_{\alpha,\beta}(\Gamma))$ and $a\in C_0(\Gamma)$, we get
$$\pi(\alpha^\#(a))T(\xi)=T(\xi)\pi(\beta^\#(a))$$
where $(T,\pi)$ is the universal representation coisometric on $J_{\cE_{Q_{\alpha,\beta}(\Gamma)}}$
\index{J$\mbox{}_\cE$} onto $\cO_{\alpha,\beta}(\Gamma).$\\ 
\index{Cuntz-Pimsner Algebra}\index{$\cO_{\alpha,\beta}(\Gamma)$}
\pf Let $a\in C_0(\Gamma)$ and $\xi\in C(\Omega_{\alpha,\beta}(\Gamma))$ then 
\begin{align*}
(\alpha^\#(a)\cdot \xi)(x,y)&=\alpha^\#(a)(y)\xi(x,y)\\
&=a(\alpha(y))\xi(x,y)\\
&=\xi(x,y)a(\beta(x))\\
&=(\xi\cdot\beta^\#(a))(x,y)
\end{align*}
for each $(x,y)\in\Omega_{\alpha,\beta}(\Gamma)$ and hence,
$$\pi(\alpha^\#(a))T(\xi)=T(\alpha^\#(a)\cdot\xi)=T(\xi\cdot\beta^\#(a))=T(\xi)\pi(\beta^\#(b)).$$\sq
\end{proposition}

\begin{theorem}\label{Gen2} 
\rm Suppose $\Gamma$ is a compact group with $\alpha,\beta\in\mbox{End}(\Gamma)$ and the left action 
$\phi$ is injective (i.e., $Q_{\alpha,\beta}(\Gamma)$ is sinkless). If there exists an orthonormal basis 
\index{Basis! Orthonormal} $\{u_i\}_{i=1}^n$ then
$\cO_{\alpha,\beta}(\Gamma)$ \index{$\cO_{\alpha,\beta}(\Gamma)$}\index{Cuntz-Pimsner Algebra}
is the universal $C^*$-algebra generated by $A=C(\Gamma)$ and isometries $\{S_i\}_{i=1}^n$
subject to
\begin{enumerate}
\item $\alpha^\#(a)S_i=S_i\beta^\#(a)$ for each $i=1,...,n$ and $a\in C(\Gamma)$
\item $S_i^*S_j=\delta_i^j$
\item $\sum_{i=1}^n S_iS_i^*=1.$
\end{enumerate}
If, in addition, there exist $\{a_i\}_{i=1}^n\subseteq C(\Gamma)$ and $u\in\{u_i\}_{i=1}^n$ such that 
$a_i\cdot u = u_i$ for each $i=1,...,n$
then $\cO_{\alpha,\beta}(\Gamma)$ is the universal $C^*$-algebra generated by $A=C(\Gamma)$ and an isometry $S$
subject to
\begin{enumerate}
\item $\alpha^\#(a)S=S\beta^\#(a)$ for each $i=1,..,n$ and $a\in C(\Gamma)$
\item $S^*a_iS=\delta_0^i$
\item $\sum_{i=1}^n a_iSS^*a_i^*=1.$
\end{enumerate}
\pf This is none other than a special case of Theorem \ref{BasisGen}.\\\sq
\end{theorem}

\begin{corollary}\rm Let $\Gamma$ be a compact group with $\alpha,\beta\in\mbox{End}(\Gamma)$ and $\phi$ injective. 
If there exists $u\in\cE_{\alpha,\beta}(\Gamma)$\index{$\cE_{\alpha,\beta}(\Gamma)$}  
where $u\innprod{u}{\xi}=\xi$ for each $\xi\in\cE_{\alpha,\beta}(\Gamma)$
(i.e., $\theta_{u,u}=1$) then $\cO_{\alpha,\beta}(\Gamma)$ is the universal $C^*$-algebra generated by $A=C(\Gamma)$ and 
unitary $S$ such that
$$\alpha^\#(a)S=S\beta^\#(a)$$ 
for all $a\in C(\Gamma).$
\end{corollary}

\begin{remark}\rm This corollary suggests a similarity to crossed products by $\zZ$. In some sense, what we have here is a
``crossed product'' by two endomorphisms. In fact, none of these last results are exclusive to topological group relations. Given
a second countable locally compact Hausdorff space $X$, define a topological quiver 
$Q=(X,\Omega_{\alpha,\beta}(X),r,s,\lambda)$ with
$$\Omega_{\alpha,\beta}(X)=\{(x,y)\st \alpha(y)=\beta(x)\}$$
for continuous maps $\alpha,\beta$ on $X$ and appropriate measures $\lambda.$ Then similar conclusions remain true. 
\end{remark}
\section{Some $C^*$-subalgebras of $\cO_{F,G}(\zT^d)$}By Proposition \ref{equiv}, we may assume (for the remainder of this paper) that $F=\diag{a_1,...,a_d}$ with $a_j>0.$ Let $G\in M_d(\zZ)$ with $j$-th row denoted $G_j$ and $\det G\ne0$, then recall the definitions of $\{U_j\}_{j=1}^d$ and $S$ from Section 3.3. $U_j\in C(\zT^d)$ is the unitary defined by $U_j(x)=x_j$ for each $x=(x_k)_{k=1}^d\in\zT^d$, $j\in\{1,...,d\}$ and 
$S$ is the isometry defined as the image of $u_0=1\in C(\Omega_{F,G}(\zT^d))$ into $\cO_{F,G}(\zT^d).$ By Theorem \ref{TdQuiverGen}, the $C^*$-algebra $\cO_{F,G}(\zT^d)$ is the universal $C^*$-algebra generated by the commuting 
unitaries $\{U_j\}_{j=1}^d$ and isometry $S$ satisfying:
\begin{enumerate}
\item $S^*U^\nu S=\delta_0^\nu,$ for all $\nu\in\fI(F)$,
\item $U_j^{a_j}S=SU^{G_j},$ for all $j=1,...,d$ and
\item $\sum_{\nu\in\fI(F)} U^\nu SS^*U^{-\nu}=1$
\end{enumerate} 
where $\fI(F)=\{\nu=(\nu_j)_{j=1}^d\in\zZ^d\st 0\le \nu_j\le a_j-1\}.$

\begin{lemma}\label{Subalg4.1}\rm Given $k\in\zN$ and $\abs{\det G}=1$, the $C^*$-algebra 
$C^*(\{U_j\}_{j=1}^d, S^k)$ viewed as a $C^*$-subalgebra of $\cO_{F,G}(\zT^d)$ is a quotient of 
$\cO_{F^k,G^k}(\zT^d)$.\\
\pf Let $\tilde{S}_\nu=U^\nu S^k$ for $\nu\in\fI(F^k)=\{(\nu_j\}_{j=1}^d\st 0\le\nu_j\le a_j^k-1\}$ $(\tilde{S}=S^k)$ and so
 \begin{enumerate}
\item $\tilde{S}^*U^{\nu^\prime}\tilde{S}=(S^k)^*U^{\nu^\prime}S^k=\delta_0^{\nu^\prime},$ 
for each $\nu^\prime\in\fI(F^k).$
\end{enumerate}
Also, by (3) of Theorem \ref{TdQuiverGen},
$$U^{a_jG_j}S=[\prod_{m=1}^d U_m^{a_jb_{jm}}]S=S\prod_{m,l=1}^d U_l^{b_{jm}b_{ml}}= SU^{(G^2)_j}.$$ 
Thus,
\begin{enumerate}
\item[2.] $U_j^{a_j^k}\tilde{S}=SU^{a_j^{k-1}G_j}S^{k-1}=...=S^kU^{(G^k)_j}=\tilde{S}U^{(G^k)_j}$
\end{enumerate}
For $k=2$, we may write $\nu = \nu_0+\mu\cdot\nu_1$ [point-wise product and sum] 
where $\mu=(a_1,...,a_d)$, $\nu_i\in\fI(F).$ Next, note for any $\omega\in\fI(F),$ there exist unique $\nu_1\in\fI(F)$
and unique $k_{\nu_1}\in\zZ^d$ such that $(\nu_1+k_{\nu_1}\cdot\mu)G=\omega$ (matrix multiplication on the right)
by the bijectivity of $G.$ Hence, 
\begin{align*}
\sum_{\nu_1\in\fI(F)}U^{\nu_1G}SS^*U^{-\nu_1G}&
=\sum_{\nu_1}U^{(\nu_1+k_{\nu_1}\cdot\mu)G}U^{(k_{\nu_1}\cdot\mu)G}S
S^*U^{-(k_{\nu_1}\cdot\mu)G}U^{-(\nu_1+k_{\nu_1}\cdot\mu)G}\\
&=\sum_{\nu_1}U^{(\nu_1+k_{\nu_1}\cdot\mu)G}SU^{k_{\nu_1}G^2}U^{-k_{\nu_1}G^2}
S^*U^{-(\nu_1+k_{\nu_1}\cdot\mu)G}\\
&=\sum_{\nu_1}U^{(\nu_1+k_{\nu_1}\cdot\mu)G}SS^*U^{-(\nu_1+k_{\nu_1}\cdot\mu)G}\\
&=\sum_{\nu_1^\prime\in\fI(F)} U^{\nu_1^\prime G}SS^*U^{-\nu_1^\prime G}=1.
\end{align*}
This gives us  
\begin{align*}
\sum_{\nu\in\fI(F^k)} U^\nu\tilde{S}\tilde{S}^*U^{-\nu}
&=\sum_{\nu_0,\nu_1\in\fI(F)}U^{\nu_0}U^{\mu\nu_1}S^2(S^2)^*U^{-\mu\nu_1}U^{-\nu_0}\\
&=\sum_{\nu_0,\nu_1}U^{\nu_0}S[U^{\sum_{j=1}^d(\nu_1)_jG_j}SS^*U^{-\sum_{j=1}^d(\nu_1)_jG_j}]S^*U^{-\nu_0}\\
&=\sum_{\nu_0}U^{\nu_0}SS^*U^{-\nu_0}=1.
\end{align*}
A similar argument can be used to show $\sum_{\nu\in\fI(F^k)}\tilde{S}_\nu\tilde{S}_\nu^*=1$ where $k\in\zN$. Hence, by Theorem \ref{TdQuiverGen}, universality implies that there exists a surjective homomorphism of $\cO_{F^k,G^k}(\zT^d)$ onto $C^*(\{U_j\}_{j=1}^d, S^k).$\\\sq
\end{lemma}

\begin{lemma}\label{Subalg}\rm If $k_j\vert a_j$ for each $j=1,...,d$ and $\gcd(\det F,\det G)=1$ then 
$$C^*(\{U_j^{k_j}\}_{j=1}^d,S)=\cO_{F,G}(\zT^d).$$
Next assume $k_j\vert a_j$ for all $j=1,...,d$. Then there are $p_j\in\zN$ such that $a_j=k_jp_j$, hence
$$U_j^{a_j}=(U_j^{k_j})^{p_j}\in C^*(\{U_j^{k_j}\}_{j=1}^d, S)$$
for all $j=1,...,d$. Further note $U^{G_j}=S^*U_j^{a_j}S\in C^*(\{U_j^{k_j}\}_{j=1}^d, S).$

Let $Q=(q_{ij})_{ij}=[Q_1 ... Q_d]^T=\mbox{adj } G\in M_d(\zZ)$, the adjugate of $G$ and $G=[G^\prime_1 ... G^\prime_d]$ 
so that $QG=GQ=(\det G)1_d$ and hence, $Q_jG^\prime_m=(\det G)\delta_j^m$. 
Then given $t\in\zT^d$ and $j\in\{1,...,d\}$, we have
$$\sum_{l=1}^d q_{jl}G_lt=\sum_{m=1}^d\sum_{l=1}^d (q_{jl}b_{lm})t_m=\sum_m Q_jG^\prime_mt_m=(\det G)t_j$$
showing that
$$U_j^{\det G}=\prod_l (U^{G_l})^{q_{jl}}\in C^*(\{U_j^{k_j}\}_{j=1}^d, S).$$
Finally, since $\gcd(\det F,\det G)=1$, we have $\gcd(a_j, \det G)=1$ for each $j$. Hence, there exists $p,q\in\zZ$ such that $1=pa_j+q\det G$ and 
$$U_j=U_j^{pa_j+q\det G}=(U_j^{a_j})^p(U_j^{\det G})^q\in C^*(\{U_j^{k_j}\}_{j=1}^d,S)$$
for each $j=1,...,d.$
We have now shown that the generators $\{U_j\}_{j=1}^d$ and $S$ are in $C^*(\{U_j^{k_j}\}_{j=1}^d, S)$. Therefore,
$$\cO_{F,G}(\zT^d)\subseteq C^*(\{U_j^{k_j}\}_{j=1}^d,S)\subseteq \cO_{F,G}(\zT^d)$$
and hence, we have equality.\\\sq
\end{lemma}

\begin{proposition}\rm Let $\{k_j\}_{j=1}^d\subset \zN$ and $\gcd(\det F,\det G)=1.$ 
Then, as a $C^*$-subalgebra of $\cO_{F,G}(\zT^d)$, $C^*(\{U_j^{k_j}\}_{j=1}^d, S)$ is a quotient of $\cO_{F,G}(\zT^d).$\\
\pf If $\gcd(k_j,a_j)\ge 2$ then represent $k_j=p^{\alpha_j}p_j$ and $a_j=p^{\beta_j}q_j$ where $p$ is prime and $p_j, q_j$
have no factors of $p$. Suppose $\alpha_j>\beta_j$ then $k_jq_j=p^{\alpha_j-\beta_j}p_ja_j$ and
$$U^{p^{\alpha_j-\beta_j}p_jG_j}=S^*SU^{p^{\alpha_j-\beta_j}p_jG_j}=S^*U_j^{k_jq_j}S\in C^*(\{U_j^{k_j}\}_j,S).$$
As in Lemma \ref{Subalg}, we obtain
$$U_j^{p^{\alpha_j-\beta_j}p_j\det G}\in C^*(\{U_j^{k_j}\}_j, S)$$
and using $\gcd(\det F,\det G)=1$,
$$U_j^{p^{\alpha_j-\beta_j}p_j}\in C^*(\{U_j^{k_j}\}_j,S).$$

Suppose $\alpha_j\le\beta_j$ then $k_jp^{\beta_j-\alpha_j}q_j=p_ja_j,$ so similarily
$$U_j^{p_j}\in C^*(\{U_j^{k_j}\}_j,S).$$
These cases show us that we may assume $\gcd(k_j, a_j)=1$ for all $j=1,...,d.$

Suppose $\gcd(k_j,a_j)=1$ for all $j=1,...,d$, then 
$$\{(\lambda_1k_1 \mod a_1,...,\lambda_dk_d \mod a_d)\st \lambda_j\in\{0,...,a_j-1\}\}=\oplus_{j=1}^d\zZ_{a_j}.$$
Hence, given $\lambda\in\fI(F),$ there exists $l_\lambda\in\fI(F)$ and $p_\lambda\in\zZ^d$ such that
$$\lambda k =l_\lambda +\mu p_\lambda\qquad(\mbox{that is, }\lambda_jk_j=l_j+a_jp_j)$$
where $\mu=(a_1,...,a_d)$.
Set $$\tilde{S}_\lambda=U^{\lambda k}S\in C^*(\{U_j^{k_j}\}_j,S).$$ 
Then
\begin{align*}
\sum_\lambda \tilde{S}_\lambda\tilde{S}_\lambda^*&=\sum_\lambda (U^{\lambda k}S)(U^{\lambda k}S)^*\\
&=\sum_\lambda U^{l_\lambda}SS^*U^{-l_\lambda}\\
&=\sum_{l_\lambda\in\fI(F)} S_{l_\lambda}S^*_{l_\lambda}=1
\end{align*}
Furthermore, by setting $W_j$ to be the full spectrum unitary $U_j^{k_j}$, we obtain 
\begin{enumerate}
\item $\tilde{S}_\nu=W^\nu\tilde{S}$,
\item $\tilde{S}_\nu^*\tilde{S}_{\nu^\prime}=\delta_\nu^{\nu^\prime}$,
\item $W_j^{a_j}\tilde{S}=\tilde{S}W^{G_j}$, and
\item $\sum_\nu W^{\nu}\tilde{S}\tilde{S}^*W^{-\nu}=1.$
\end{enumerate}
That is, noting $S=\tilde{S}_{(0,...,0)}$, $C^*(\{U_j^{k_j}\}_{j=1}^d,S)$ is the universal $C^*$-algebra generated by commuting full spectrum unitaries
$W_j$ and isometries $\{\tilde{S}_\nu\}_{\nu\in\fI(F)}$ such that 
\begin{enumerate}
\item $\tilde{S}_\nu=W^\nu\tilde{S}$,
\item $\tilde{S}_\nu^*\tilde{S}_{\nu^\prime}=\delta_\nu^{\nu^\prime}$,
\item $W_j^{a_j}\tilde{S}=\tilde{S}W^{G_j}$, and
\item $\sum_\nu W^{\nu}\tilde{S}\tilde{S}^*W^{-\nu}=1.$
\end{enumerate}
By Theorem \ref{TdQuiverGen}, the universality of $\cO_{F,G}(\zT^d)$ implies that 
there exists a surjective homomorphism of $\cO_{F,G}(\zT^d)$ onto $C^*(\{U_j^{k_j}\}_{j=1}^d, S).$\\\sq
\end{proposition}
\section{Spatial Structure of $\cO_{F,G}(\zT^d)$}We now turn our attention to a colimit structure of $\cO_{F,G}(\zT^d).$

\begin{lemma}\label{OStruc}\rm 
Let $\cO=\overline{\mbox{span}}\{S_\alpha U^\nu S_\beta^*\st \alpha=(\alpha_j)_{j=1}^{k_1}, 
\beta=(\beta_j)_{j=1}^{k_2}, k_1,k_2\in\zN, \alpha_j,\beta_j\in\fI(F), \nu\in\zZ^d\}$ where $S_\alpha=S_{\alpha_1}\cdot\cdot\cdot S_{\alpha_{k_1}}.$ Then
$$\cO=\cO_{F,G}(\zT^d)$$
\pf Begin by noting $S_{\alpha_i}=U^{\alpha_i}S.$ With the notation $\alpha_i=(\alpha_{ij})_{j=1}^d\in\fI(F)$
$$U_jS_{\alpha_i}=U^{\alpha_i+e_j}S=
\begin{cases}
S_{\alpha_i+e_j} & \text{if $\alpha_{ij}<a_j-1$} \\
S_{\alpha_i-(a_j-1)e_j}U^{G_j}& \text{if $\alpha_{ij}=a_j-1$}
\end{cases}$$
where $e_j$ is the vector with $1$ in the $j$-th coordinate and $0$ elsewhere.
Also,\begin{enumerate}
\item $1=\sum_\nu S_\nu S^*_\nu\in\cO,$
\item $U_j=\sum_\nu S_{\nu+e_j}S^*_\nu\in\cO,$
\item $S=\sum_\nu SS_\nu S^*_\nu\in\cO,$ and
\item $S^*=\sum_\nu S_\nu S^*_\nu S^*\in\cO.$
\end{enumerate} 
Finally, 
$$U_j^*S=U_j^{a_j-1}SU^{-G_j}=S_{(a_j-1)e_j}U^{-G_j};$$
hence, any word in $\{U_j, U_j^*\}_j, S$ and $S^*$ may be expressed as the sum of elements in $\cO$.\\\sq
\end{lemma}

Given $\alpha=(\alpha_j)_{j=1}^{k}$ where $\alpha_j\in\fI(F)$\index{$\fI(F)$} 
and $k,$ the ``length of $\alpha,$'' is denoted by $\abs{\alpha}.$
Let $\gamma:\zT\to\mbox{Aut}(\cO_{F,G}(\zT^d))$ be the \emph{gauge action}\index{Gauge Action, $\gamma$} (see \cite{MT})
defined by $t\mapsto\gamma_t$ where 
$$\gamma_t(S)=tS\qquad\gamma_t(U_j)=U_j\mbox{ for each $j=1,...,d$}.$$
If $S_\alpha U^\nu S^*_\beta\in\cO_{F,G}(\zT^d)^\gamma$ 
(the \emph{fixed point algebra}\index{Fixed Point Algebra of the Gauge Action})\index{$\cO_{F,G}(\zT^d)^\gamma$} then
$$S_\alpha U^\nu S^*_\beta=\gamma_t(S_\alpha U^\nu S^*_\beta)=t^{\abs{\alpha}-\abs{\beta}}S_\alpha U^\nu S_\beta^*$$
implies $\abs{\alpha}=\abs{\beta}$ (since $S_\alpha U^\nu S^*_\beta\ne0.$) Thus,
$$\cO_{F,G}(\zT^d)^\gamma\supseteq\overline{\mbox{span}}\{S_\alpha U^\nu S^*_\beta\st \abs{\alpha}=\abs{\beta},\nu\in\zZ^d\}.$$
Furthermore, since $\zT$ is compact, averaging over $\gamma$ with respect to normalized Haar measure produces an
expectation of $\cO_{F,G}(\zT^d)$ onto $\cO_{F,G}(\zT^d)^\gamma.$ Let $E:\cO_{F,G}(\zT^d)\to\cO_{F,G}(\zT^d)^\gamma$ 
be the expectation onto the fixed point algebra $\cO_{F,G}(\zT^d)^\gamma$ defined by
$$E(x)=\int_\zT\gamma_t(x)\,dt$$\index{Expectation onto the Fixed Point Algebra}
for each $x\in\cO_{F,G}(\zT^d)$ (see \cite[Section 4.1]{MT}).
Then for $S_\alpha U^\nu S_\beta^*\in\cO_{F,G}(\zT^d)$ where $\alpha\in\fI(F)^{k_1},$ $\beta\in\fI(F)^{k_2}$ and $\nu\in\zZ^d$ as in Lemma 5.1, note 
\begin{align*}
E(S_\alpha U^\nu S_\beta^*)&=\int_\zT t^{\abs{\alpha}-\abs{\beta}}S_\alpha U^\nu S_\beta\, dt\\
&=\frac{S_\alpha U^\nu S_\beta^*}{2\pi}\int_{0}^{2\pi}e^{2\pi i [\abs{\alpha}-\abs{\beta}]z}dz\\
&=\begin{cases}
\frac{S_\alpha U^\nu S_\beta^*}{2\pi}[\frac{e^{2\pi i(\abs{\alpha}-\abs{\beta})z}}{2\pi i(\abs{\alpha}-\abs{\beta})}]_0^{2\pi}
&\mbox{if $\abs{\alpha}\ne\abs{\beta}$}\\
\frac{S_\alpha U^\nu S_\beta^*}{2\pi}[z]_0^{2\pi}&\mbox{if $\abs{\alpha}=\abs{\beta}$}
\end{cases}\\
&=\delta_{\abs{\alpha}}^{\abs{\beta}}S_\alpha U^\nu S_\beta^*
\end{align*}
and by Lemma 5.1, 
$$\cO_{F,G}(\zT^d)^\gamma=E(\cO_{F,G}(\zT^d))=\overline{\mbox{span}}\{S_\alpha U^\nu S^*_\beta\st \abs{\alpha}=\abs{\beta},\nu\in\zZ^d\}.$$

\begin{proposition}\label{FixedStruc}\rm Let
$\cA_k=\overline{\mbox{span}}\{S_\alpha U^\nu S^*_\beta\st \abs{\alpha}=\abs{\beta}=k, \nu\in\zZ^d\}$ 
for $k\in\zN\cup\{0\}.$ Then
$\cA_k$ is a unital $*$-subalgebra of $\cO_{F,G}(\zT^d)^\gamma$ with $\cA_k\subset \cA_{k+1}.$ Moreover,
$$\cO_{F,G}(\zT^d)^\gamma=\overline{\cup_k\cA_k}.$$
\pf Note $S^*_\alpha S_\beta=\delta_\alpha^\beta$, thus $\cA_k$ is closed under multiplication. Furthermore,
\begin{align*}
S_\alpha U^\nu S^*_\beta&=S_\alpha\sum_\mu S_\mu S^*_\mu U^\nu S_\beta^*\\
&=\sum_\mu S_{(\alpha,\mu)}U^{\nu_\mu}S^*_{(\beta, \varpi_\mu)}\in\cA_{k+1}
\end{align*} 
where $(\alpha,\mu)=(\alpha_1,...,\alpha_k, \mu)$ and $\varpi_\mu$ and $\nu_\mu$ are the appropriate vectors with
$$S^*_\mu U^\nu=S^*U^{\nu-\mu}=U^{\nu_\mu}S^*U^{-\varpi_\mu}=U^{\nu_\mu}S^*_{\varpi_\mu}.$$
Finally, $1=\sum_\nu S_\nu S^*_\nu\in\cA_1\subset\cA_k.$\\\sq
\end{proposition}

\begin{remark}\rm This result comes as little surprise by Proposition 5.7 of \cite{K} which states that, in general, the fixed point algebra coincides with an inductive limit. This aside, our result leads us to something more interesting.  
\end{remark}

\begin{lemma}\rm Let $E_{\alpha\beta}=S_\alpha S^*_\beta$ where $\alpha,\beta\in\fI(F)^k$. Then $\{E_{\alpha\beta}\}_{\alpha,\beta}$ is a system of matrix units for $M_{N^k}(\zC)$ in $\cA_k.$\\
\pf We need only check the properties of matrix units:
\begin{align*}
(1)E_{\alpha\beta}E_{\alpha^\prime\beta^\prime}&=S_\alpha S^*_\beta S_{\alpha^\prime}S^*_{\beta^\prime}
=\delta_{\beta}^{\alpha^\prime}S_\alpha S^*_{\beta^\prime}=\delta_{\beta}^{\alpha^\prime}E_{\alpha\beta^\prime},\\
(2)\sum_\alpha E_{\alpha\alpha}&=\sum_\alpha S_\alpha S_\alpha^* \\
&=\sum_{\alpha_1,...,\alpha_k}S_{\alpha_1}\cdot\cdot\cdot S_{\alpha_{k-1}}S_{\alpha_k}S_{\alpha_k}^*S_{\alpha_{k-1}}^* \cdot\cdot\cdot S_{\alpha_1}^*\\
&=\sum_{\alpha_1,...,\alpha_{k-1}}S_{\alpha_1}\cdot\cdot\cdot S_{\alpha_{k-1}}\big[\sum_{\alpha_k} S_{\alpha_k}S_{\alpha_k}^*\big]S_{\alpha_{k-1}}^* \cdot\cdot\cdot S_{\alpha_1}^*\\
&=\sum_{\alpha_1,...,\alpha_{k-1}}S_{\alpha_1}\cdot\cdot\cdot S_{\alpha_{k-1}}[1]S_{\alpha_{k-1}}^* \cdot\cdot\cdot S_{\alpha_1}^*\\
&=...\\
&=\sum_{\alpha_1}S_{\alpha_1}S_{\alpha_1}^*=1.
\end{align*}\sq
\end{lemma}

\begin{lemma}\cite[Lemma 11.4.8]{KR2} \rm If $\mathcal M$ is a type I factor and $B$ is a $C^*$-subalgebra of $\mathcal M^\prime,$ the commutant of $\mathcal M$, and $A$ is the $C^*$-algebra generated by $\mathcal M\cup B$, then there is
an isomorphism $\varphi:\mathcal M\otimes B\to A$ such that $\varphi(m\otimes b)=mb$ for elementary tensors 
$m\otimes b\in\mathcal M\otimes B.$
\end{lemma}

Note $\mathcal M_k:=\mbox{span}\{S_\alpha S_\beta^*\st \alpha,\beta\in\fI(F)^k\}\cong M_{N^k}(\zC)$ (by Lemma 5.4) is a type I factor. Let $B_k=C^*(\{U_j^{a_j^k}\}_{j=1}^d)\cong C(\zT^d)$ (since $U_j^{a_j^k}$ is a full spectrum unitary (see Remark \ref{FullSpec})), then $\cA_k$ is the $C^*$-algebra generated by $\mathcal M_k$ and $B_k.$ Indeed, note
$$U_j^{a_j^k}S_\alpha S_\beta^*=S_\alpha U^{(G^k)_j} S_\beta^*$$
for each $\alpha,\beta\in\fI(F)^k$ and $j=1,...,d$. 
Note, in the proof of Lemma 4.1, there exists $\{q_{jl}\}_{j,l=1}^d\subset\zZ$ such that 
$$\sum_{l=1}^dq_{jl}(G^k)_lt=(\det G^k)t_j$$
for each $j=1,...,d.$ Hence, $U_j^{\det G^k}\in C^*(\{U^{(G^k)_j}\}_{j=1}^d).$ 
Further note $\{U_j^{\det G^k}\}_{j=1}^d$ generates the faithful image of $C(\zT^d)$ in $\cO_{F,G}(\zT^d)$
(since $\det G^k=(\det G)^k\ne 0$ and the left action is injective.)
Thus, $\cA_k$ is the $C^*$-algebra generated by $\mathcal M_k$ and $B_k.$ 
Furthermore,
$$U_j^{a_j^k}S_\alpha S_\beta^*=S_\alpha U^{(G^k)_j}S_\beta^*=S_\alpha (S_\beta U^{-(G^k)_j})^*=S_\alpha (U_j^{-a_j^k}S_\beta)^*=S_\alpha S_\beta^*U_j^{a_j^k}$$
for each $j=1,...,d$ and $\alpha,\beta\in\fI(F)^k;$ hence, the generators of $B$ commute with the generators of $\mathcal M.$
Thus, $B$ may be viewed as a $C^*$-subalgebra of $\mathcal M^\prime.$ Hence, we use the previous lemma. 
 
\begin{theorem}\label{Tensor}\rm $$\cA_k\cong C(\zT^d)\otimes M_{N^k}(\zC)$$
\pf By Lemma 5.5, we obtain an isomorphism $\psi_k:B_k\otimes\mathcal M_k\to\cA_k.$ Let $\psi_{B_k}:B_k\to C(\zT^d)$ and $\psi_{\mathcal M_k}:\mathcal M_k\to M_{N^k}(\zC)$ be isomorphisms, then 
$$(\psi_{B_k}\otimes\psi_{\mathcal M_k})\circ(\psi_k)^{-1}:\cA_k\to C(\zT^d)\otimes M_{N^k}(\zC)$$ 
is an isomorphism.\\\sq
\end{theorem}

Hence, for the maps $\varphi_k: C(\zT^d)\otimes M_{N^k}(\zC)\cong\cA_k\hookrightarrow\cA_{k+1}\cong C(\zT^d)\otimes M_{N^{k+1}}(\zC),$ we have
$$\cO_{F,G}(\zT^d)^\gamma=\mbox{colim} (M_{N^k}(C(\zT^d)), \varphi_k).$$
We use the identification $C(\zT^d)\otimes M_{N^k}(\zC)\cong M_{N^k}(C(\zT^d)).$ 

We now define the notion of a crossed product by an endomorphism as discussed in \cite{S}.

\begin{definition}\cite{S}\label{CP} \rm A \emph{crossed product}\index{Crossed Product by an Endomorphism}
(of multiplicity 1) for an endomorphism $\rho$ of a unital $C^*$-algebra $A$ with 
$$A_\infty=\mbox{colim}_{k\in\zN} (A,\rho)\ne 0$$
is the (unique) unital $C^*$-algebra $A\rtimes_\rho\zN$ together with a homomorphism $\iota:A\to A\rtimes_\rho\zN$ with $\iota(1_A)=1_{A\rtimes\rho\zN}$, and an isometry $s\in A\rtimes_\rho\zN$ such that
\begin{enumerate}
\item $\iota(\rho(a))=s\iota(a)s^*$ for all $a\in A$
\item for every representation $\pi$ of $A$ and isometry $t$ satisfying $\pi(\rho(a))=t\pi(a)t^*$ for all $a\in A$, there exists
a non-degenerate representation $\pi\times T$ of $A\rtimes_\rho\zN$ with
$$(\pi\times T)\circ\iota=\pi\qquad\mbox{and}\qquad (\pi\times T)(s)=t.$$
\item $A\rtimes_\rho\zN$ is generated by $\iota(A)$ and $s.$ 
\end{enumerate}
\end{definition}

Let $\rho_k:C(\zT^d)\otimes M_{N^k}(\zC)\to C(\zT^d)\otimes M_{N^{k+1}}(\zC)$ be defined on basic tensors of $C(\zT^d)\otimes M_{N^k}(\zC)$ by 
$$f\otimes E_{\alpha\beta}\mapsto f\otimes E_{(0,\alpha)(0,\beta)}$$ 
where $\alpha,\beta\in\fI(F)^k$ and $f\in C(\zT^d).$ Note that since $\{E_{\alpha\beta}\}_{\alpha,\beta\in\fI(F)^k}$ is linearly independent, the definition of $\rho_k$ extends, by linearity, to all of $C(\zT^d)\otimes M_{N^k}(\zC).$

\begin{lemma}\rm Let $\rho_k:C(\zT^d)\otimes M_{N^k}(\zC)\to C(\zT^d)\otimes M_{N^{k+1}}(\zC)$ be defined by 
$$f\otimes E_{\alpha\beta}\mapsto f\otimes E_{(0,\alpha)(0,\beta)}$$
on the basic tensors, $f\otimes E_{\alpha\beta},$ of $C(\zT^d)\otimes M_{N^k}(\zC)$ where $\alpha,\beta\in\fI(F)^k$ and
$f\in C(\zT^d).$ Then $\rho_k$ is a (continuous) $*$-homomorphism for each $k\in\zN$ and the following diagram commutes
$$\begindc{\commdiag}[8]
\obj(0,6)[B]{$\cdots$}
\obj(15,6)[C]{$\mathcal C(\zT^d)\otimes M_{N^k}(\zC)$}
\obj(30,6)[D]{$\mathcal C(\zT^d)\otimes M_{N^{k+1}}(\zC)$}
\obj(45,6)[E]{$\cdots$}
\obj(0,0)[G]{$\cdots$}
\obj(15,0)[H]{$\mathcal C(\zT^d)\otimes M_{N^k}(\zC)$}
\obj(30,0)[I]{$\mathcal C(\zT^d)\otimes M_{N^{k+1}}(\zC)$}
\obj(45,0)[J]{$\cdots$}
\mor{B}{C}{$\varphi_{k-1}$}
\mor{C}{D}{$\varphi_{k}$}
\mor{D}{E}{$\varphi_{k+1}$}
\mor{G}{H}{$\varphi_{k-1}$}
\mor{H}{I}{$\varphi_{k}$}
\mor{I}{J}{$\varphi_{k+1}$}
\mor{G}{C}{$\rho_{k-1}$}
\mor{H}{D}{$\rho_{k}$}
\mor{I}{E}{$\rho_{k+1}$}
\enddc$$
\pf Note for any $\sum_{\alpha,\beta\in\fI(F)^k}f_{\alpha,\beta}\otimes E_{\alpha\beta}$ and $\sum_{\alpha^\prime,\beta^\prime\in\fI(F)^k}f_{\alpha^\prime,\beta^\prime}\otimes E_{\alpha^\prime\beta^\prime}$
\begin{align*}
\rho_k((\sum_{\alpha,\beta\in\fI(F)^k}f_{\alpha,\beta}\otimes E_{\alpha\beta})&(\sum_{\alpha^\prime,\beta^\prime\in\fI(F)^k}f_{\alpha^\prime,\beta^\prime}\otimes E_{\alpha^\prime\beta^\prime}))\\
&=\rho_k(\sum_{\alpha,\beta,\alpha^\prime,\beta^\prime\in\fI(F)^k}f_{\alpha,\beta}f_{\alpha^\prime,\beta^\prime}\otimes E_{\alpha\beta}E_{\alpha^\prime\beta^\prime})\\
&=\rho_k(\sum_{\alpha,\beta^\prime\in\fI(F)^k}f_{\alpha,\beta}f_{\beta,\beta^\prime}\otimes E_{\alpha\beta^\prime})\\
&=\sum_{\alpha,\beta^\prime\in\fI(F)^k}f_{\alpha,\beta}f_{\beta,\beta^\prime}\otimes E_{(0,\alpha)(0,\beta^\prime)})\\
&=(\sum_{\alpha,\beta\in\fI(F)^k}f_{\alpha,\beta}\otimes E_{(0,\alpha)(0,\beta)})(\sum_{\alpha^\prime,\beta^\prime\in\fI(F)^k}f_{\alpha^\prime,\beta^\prime}\otimes E_{(0,\alpha^\prime)(0,\beta^\prime)})\\
&=\rho_k(\sum_{\alpha,\beta\in\fI(F)^k}f_{\alpha,\beta}\otimes E_{\alpha\beta})\rho_k(\sum_{\alpha^\prime,\beta^\prime\in\fI(F)^k}f_{\alpha^\prime,\beta^\prime}\otimes E_{\alpha^\prime\beta^\prime})\\
\end{align*}

Recall that $\mathcal A_k\subset \mathcal A_{k+1}.$ Let $i_k$ denote the inclusion and recall, 
for any $S_\alpha U^\nu S_\beta^*\in\mathcal A_k$
$$i_k(S_\alpha U^\nu S_\beta^*)=S_\alpha(\sum_{\mu\in\fI(F)} S_\mu S_\mu^*)U^\nu S_\beta^*
=\sum_{\mu\in\fI(F)} S_{(\alpha,\mu)}U^{\nu_\mu}S_{(\beta,\varpi_\mu)}^*$$
for appropriate $\nu_\mu\in\zZ^d$ and $\varpi_\mu\in\fI(F)$ as in Proposition \ref{FixedStruc}. Thus
for $\varphi_k$ defined in the paragraph proceeding Theorem \ref{Tensor}, 
$$\varphi_k(U^\nu\otimes E_{\alpha\beta})=\sum_{\mu\in\fI(F)} U^{\nu_\mu}\otimes E_{(\alpha,\mu)(\beta,\varpi_\mu)}$$
for the appropriate $\nu_\mu\in\zZ^d$ and $\varpi_\mu\in\fI(F).$
Furthermore, since 
\begin{align*}
i_{k+1}(S_{(0,\alpha)} U^\nu S_{(0,\beta)}^*)&=SS_\alpha(\sum_{\mu\in\fI(F)} S_\mu S_\mu^*)U^\nu S_\beta^*S^*\\
&=S(\sum_{\mu\in\fI(F)} S_{(\alpha,\mu)}U^{\nu_\mu}S_{(\beta,\varpi_\mu)}^*)S^*\\
&=\sum_{\mu\in\fI(F)} S_{(0,\alpha,\mu)}U^{\nu_\mu}S_{(0,\beta,\varpi_\mu)}^*
\end{align*}
we obtain
$$\varphi_{k+1}(U^\nu\otimes E_{(0,\alpha)(0,\beta)})
=\sum_{\mu\in\fI(F)} U^{\nu_\mu}\otimes E_{(0,\alpha,\mu)(0,\beta,\varpi_\mu)}.$$
Hence,
\begin{align*}
(\rho_{k+1}\circ\varphi_k)(U^\nu\otimes E_{\alpha\beta})
&=\rho_{k+1}(\sum_{\mu\in\fI(F)} U^{\nu_\mu}\otimes E_{(\alpha,\mu)(\beta,\varpi_\mu)})\\
&=\sum_{\mu\in\fI(F)} U^{\nu_\mu}\otimes E_{(0,\alpha,\mu)(0,\beta,\varpi_\mu)}\\
\end{align*}
and
\begin{align*}
(\varphi_{k+1}\circ\rho_k)(U^\nu\otimes E_{\alpha\beta})&=\varphi_{k+1}(U^\nu\otimes E_{(0,\alpha)(0,\beta)})\\
&=\sum_{\mu\in\fI(F)} U^{\nu_\mu}\otimes E_{(0,\alpha,\mu)(0,\beta,\varpi_\mu)}\\
&=(\rho_{k+1}\circ\varphi_k)(U^\nu\otimes E_{\alpha\beta})
\end{align*}
for all $U^\nu \otimes E_{\alpha\beta}\in C(\zT^d)\otimes M_{N^k}(\zC).$ The result now follows from the linearity and continuity of $\rho_k$ and $\varphi_k.$ \\\sq
\end{lemma}

Let $j_k:\mathcal A_k\hookrightarrow \cO_{F,G}(\zT^d)^\gamma$ be the natural inclusion $*$-homomorphism for each $k\in\zN$ and  let $\tilde\rho_k:\mathcal A_k\to \cO_{F,G}(\zT^d)^\gamma$ be the composition $j_{k+1}\circ\Psi_{k+1}^{-1}\circ\rho_k\circ\Psi_k$ for $k\in\zN$
where $\Psi_k$ is the isomorphism in Theorem \ref{Tensor}. 
Then the solid diagram
$$\begindc{\commdiag}[3]
\obj(40,25)[A]{$\cO_{F,G}(\zT^d)^\gamma$}
\obj(0,0)[F]{$\cdots$}
\obj(20,0)[G]{$\mathcal A_{k-1}$}
\obj(40,0)[H]{$\mathcal A_k$}
\obj(60,0)[I]{$\mathcal A_{k+1}$}
\obj(80,0)[J]{$$}
\obj(90,0)[L]{$\cdots$}
\obj(110,0)[K]{$\cO_{F,G}(\zT^d)^\gamma$}
\mor{G}{A}{$\tilde\rho_{k-1}$}
\mor{H}{A}{$\tilde\rho_{k}$}[-1,0]
\mor{I}{A}{$\tilde\rho_{k+1}$}[-1,0]
\mor{K}{A}{$\rho$}[-1,1]
\mor{F}{G}{$i_{k-2}$}[-1,0]
\mor{G}{H}{$i_{k-1}$}[-1,0]
\mor{H}{I}{$i_{k}$}[-1,0]
\mor{I}{J}{$i_{k+1}$}[-1,0]
\mor{L}{K}{$$}
\enddc$$
commutes and hence, there exists a (unique) $*$-endomorphism $\rho$ on $\cO_{F,G}(\zT^d)^\gamma$ 
making the entire diagram commute.  A diagram chase ensures that
$$\rho(S_\alpha f S_\beta^*)=S_{(0,\alpha)}fS_{(0,\beta)}^*$$
for $S_\alpha f S_\beta^*\in\cO_{F,G}(\zT^d)^\gamma.$

\begin{theorem}\label{CrossProduct}\rm \index{$\cO_{F,G}(\zT^d)$! Crossed Product by an Endomorphism} 
$$\cO_{F,G}(\zT^d)\cong\cO_{F,G}(\zT^d)^\gamma\rtimes_\rho\zN$$
where $\rho\in\mbox{End}(\cO_{F,G}(\zT^d)^\gamma)$ is defined above.\\
\pf Let $(\iota,s):\cO_{F,G}(\zT^d)^\gamma\to\cO_{F,G}(\zT^d)^\gamma\rtimes_\rho\zN$ 
be the universal representation as described in Definition \ref{CP}. Then let
$\iota_0:\cO_{F,G}(\zT^d)^\gamma\to\cO_{F,G}(\zT^d)$ be the natural inclusion and let $S\in\cO_{F,G}(\zT^d)$ denote
the isometry$S_{(0,0,...,0)}$ defined in Theorem \ref{TdQuiverGen}. Note
$$\iota_0(\rho(S_\alpha U^\nu S_\beta^*))=S_{(0,\alpha)}U^\nu S_{(0,\beta)}^*=SS_\alpha U^\nu S_\beta^*S^*=S\iota_0(S_\alpha U^\nu S_\beta^*)S^*$$
for all generators $S_\alpha U^\nu S_\beta^*\in\cO_{F,G}(\zT^d)^\gamma.$ Hence,
$$\iota_0(\rho(x))=S\iota_0(x)S^*$$
for each $x\in\cO_{F,G}(\zT^d)^\gamma.$ 
Furthermore, note $C^*(\cO_{F,G}(\zT^d)^\gamma,S)=\cO_{F,G}(\zT^d).$ Therefore,
there exists a homomorphism $$\tau:\cO_{F,G}(\zT^d)^\gamma\rtimes_\rho\zN\to\cO_{F,G}(\zT^d)$$ such that
$$\tau\circ\iota=\iota_0\qquad\mbox{and}\qquad\tau(s)=S.$$

Recall by Lemma \ref{OStruc}, there exists $x\in\cO_{F,G}(\zT^d)^\gamma$ such that $\rho(x^\nu)=U^\nu$ for each $\nu\in\fI(F)$.\index{$\fI(F)$} Then
\begin{enumerate}
\item Given $j\in\{1,...,d\},$ 
$$\iota(U_j)^{a_j}s=\iota(\rho(x))^{a_j}s=s\iota(x)^{a_j}s^*s=s\iota(x^{a_j})$$ and
$$x^{a_j}=S^*Sx^{a_j}S^*S=S^*\rho(x)^{a_j}S=S^*U^{a_j}S=U^{G_j}.$$
Hence, $$\iota(U_j)^{a_j}s=s\iota(U)^{G_j}$$ for each $j.$
\item Given $\nu\in\fI(F)$,
$$s^*\iota(U)^\nu s=s^*\iota(\rho(x))^\nu s=s^*s\iota(x)^\nu s^*s=\iota(x)^\nu$$ and
$$x^\nu=S^*Sx^\nu S^*S=S^*\rho(x)^\nu S=S^*U^\nu S=\delta_0^\nu.$$
Hence, $$s^*\iota(U)^\nu s=\delta_0^\nu$$ for each $\nu\in\fI(F).$
\item Finally, $$\sum_\nu \iota(U)^\nu ss^*\iota(U)^{-\nu}=\sum_\nu\iota(U)^\nu\iota(\rho(1))\iota(U)^{-\nu}
=\sum_\nu \iota(U^\nu SS^*U^{-\nu})=\iota(1)=1.$$
\end{enumerate} Therefore, by Theorem \ref{TdQuiverGen}, there exists a homomorphism 
$$\sigma:\cO_{F,G}(\zT^d)\to\cO_{F,G}(\zT^d)^\gamma\rtimes_\rho\zN$$
such that 
$$U_j\mapsto \iota(U_j)\mbox{ for each $j=1,...,d$}\qquad\mbox{and}\qquad S\mapsto s.$$ 
Then
$$(\tau\circ\sigma)(U_j)=(\tau\circ\iota)(U_j)=\iota_0(U_j)=U_j$$
for each $j=1,...,d$ and
$$(\tau\circ\sigma)(S)=\tau(s)=S.$$
Hence, $\tau\circ\sigma=\mbox{id}_{\cO_{F,G}(\zT^d)}$. Furthermore, for each 
$S_\alpha U^\nu S_\beta^*\in\cO_{F,G}(\zT^d)^\gamma,$
\begin{align*}
(\sigma\circ\tau)&(\iota(S_\alpha U^\nu S_\beta^*))=\sigma(S_\alpha U^\nu S_\beta^*)\\
&=\sigma(U^{\alpha_1}SU^{\alpha_2}S\cdots U^{\alpha_k}SU^\nu S^*U^{-\beta_k}S^*U^{-\beta_{k-1}}\cdots S^*U^{-\beta_1}\\
&=\iota(U^{\alpha_1})s\iota(U^{\alpha_2})s\cdots \iota(U^{\alpha_k})s\iota(U^\nu)s^*\iota(U^{-\beta_k})s^*\iota(U^{-\beta_{k-1}})\cdots s^*\iota(U^{-\beta_1})\\
&=\iota(U^{\alpha_1})s\iota(U^{\alpha_2})s\cdots\iota(U^{\alpha_{k-1}})S^*\iota(U^{\alpha_k}SU^\nu S^*U^{-\beta_k}))s^*\iota(U^{-\beta_{k-1}})\cdots s^*\iota(U^{-\beta_1})\\
&=\cdots=\iota(S_\alpha U^\nu S_\beta^*)
\end{align*}
and $$(\sigma\circ\tau)(s)=\sigma(S)=s.$$ Thus, by Lemma \ref{OStruc}, 
$\sigma\circ\tau=\mbox{id}_{\cO_{F,G}(\zT^d)^\gamma\rtimes_\rho\zN}$ and so,
$$\cO_{F,G}(\zT^d)\cong \cO_{F,G}(\zT^d)^\gamma\rtimes_\rho\zN\cong(\mbox{colim}_k (M_{N^k}(C(\zT^d)))\rtimes_\rho\zN$$\sq
\end{theorem}

\section{Acknowledgements}This paper is the first product of my doctoral dissertation at the University of Calgary. It
would not have been possible without the guidance and help of Dr. Berndt Brenken and Dr. Marcelo Laca for their 
insight on corrections and enhancements of the material discussed here.

I am also indebted to NSERC, the Department of Mathematics and Statistics at the
University of Calgary and the Department of Mathematics and Statistics at the University of Regina 
in providing funding to finance my mathematical studies.



\begin{thebibliography}{50}
\bibitem{aHR} an Huef, A., and Raeburn, I., \emph{The ideal structure of Cuntz-Krieger algebras}, Ergodic
Theory Dyn. Sys. 17 (1997) 611-624.

\bibitem{BHRS} Bates, T., Hong, J. H., Raeburn, I., and Szymanski, W., \emph{The ideal structure of $C*$-algebras of
infinite graphs}, Illinois J. Math 46 (2002), 1159-1176.

\bibitem{BPRS} Bates, T., Pask, D., Raeburn, I., and Szymanski, W., \emph{The $C^*$-algebras of row-finite graphs},
New York J. Math. 6 (2000), 307-324.

\bibitem{B1} Blackadar, B., K-theory for Operator Algebras, M. S. R. I. Monographs, vol. 5, Springer-Verlag, Berlin and New York, 1986.

\bibitem{B2} Blackadar, B., Operator Algebras: Theory of C*-Algebras and von Neumann Algebras, Encyclopaedia of Mathematical Sciences, vol. 122, Springer-Verlag, Berlin, 2006.

\bibitem{BB1} Brenken, B., \emph{The local product structure of expansive automorphisms of solenoids and their associated $C^*$-algebras}, Can. J. Math., 48, (1996), 692-709.

\bibitem{BB2} Brenken, B., \emph{$C^*$-algebras associated with topological relations}, J. Ramanujan Math. Soc. 19, No. 1 (2004), 1-21.

\bibitem{BB2.5} Brenken, B., \emph{Endomorphisms of type I von Neumann algebras with discrete center}, J. Operator Theory 51 (2004), no. 1, 19-34.

\bibitem{BB3} Brenken, B., \emph{The isolated ideal of a correspondence associated with a topological quiver}. New York J. Math., 12, (2006), 1-16.

\bibitem{BB4} Brenken, B., \emph{A Dynamical Core for Topological Directed Graphs}, Munster J. of Math. 3 (2010), 111-144. 

\bibitem{BB5} Brenken, B., \emph{Topological Quivers as Multiplicity Free Relations}, Math. Scand., 106, (2010), 217-242.

\bibitem{BO} Brown, N.P., and Ozawa, N., $C^*$-Algebras and Finite-Dimensional Approximations, Graduate Studies in Mathematics, 88. Amer. Math. Soc. Providence, Rhode Island, 2008.

\bibitem{C} Conway, J.B., A Course in Functional Analysis, Second Edition, Graduate Texts in Mathematics, 96. Springer, New York, 1990. 

\bibitem{Cu} Cuntz, J., \emph{Simple $C^*$-algebras generated by isometries}. Comm. Math. Phys. 57 (1977), no. 2, 173-185. 

\bibitem{CK} Cuntz, J., and Krieger, W., \emph{A class of $C^*$-algebras and topological Markov chains}, Inventiones Math., 56 (1980), 251-268.

\bibitem{D} Davidson, K., $C^*$-Algebras by Example, Fields Institute Monograph, Amer. Math. Soc. Providence, Rhode Island, 1996.

\bibitem{Dea} Deaconu, V., \emph{Groupoids associated with endomorphisms}, Trans. Amer. Math. Soc. 347 (1995),
1779-1786.

\bibitem{Dea2} Deaconu, V., \emph{A path model for circle algebras}, J. Operator Theory 34 (1995), 57-89.

\bibitem{Dea3} Deaconu, V., \emph{Generalized Cuntz-Krieger algebras}, Proc. Amer. Math. Soc. 124 (1996), 3427-3435.

\bibitem{Dea4} Deaconu, V., \emph{Continuous graphs and $C^*$-algebras}, in Operator theoretical methods (Timi¸soara,
1998), 137-149, Theta Found., Bucharest, 2000.

\bibitem{DM} Deaconu, V., and Muhly, P., \emph{$C^*$-algebras associated with branched coverings},
Proc. Amer. Math. Soc. 129 (2001), 1077-1086.

\bibitem{DF} Dummit, D.S., and Foote, R.M., Abstract Algebra, 3rd edition, Wiley and Sons, 2004.

\bibitem{Ev} Evans, D.E., \emph{The C*-algebras of topological Markov chains}, Lecture Notes, Tokyo Metropolitan
University (1983).

\bibitem{EaHR} Exel, R., an Huef, A., and Raeburn, I., \emph{Purely Infinite Simple $C^*$-algebras associated to Integer Dilation Matrices},  Indiana Univ. Math. J. 60 (2011), no. 3, 1033-1058. 

\bibitem{FLR} Fowler, N.J., Laca, M., and Raeburn, I., \emph{The $C^*$-algebras of infinite graphs}, Proc. Amer. Math. Soc. 8 (2000), 2319-2327.

\bibitem{FMR} Fowler, N.J., Muhly, P.S., and Raeburn, I., \emph{Representations of Cuntz-Pimsner Algerbas}, Indiana Univ. Math. J.,
52(3) (2003), 569-605.

\bibitem{FR} Fowler, N.J., Raeburn, I., \emph{The Toeplitz algebra of a Hilbert bimodule}, Indiana Univ. Math. J. 48 (1999), 155-181.

\bibitem{G} Gabriel, P., \emph{Unzerlegbare Darstellungen I} (Oberwolfach 1970), Manuscr. Math. 6 (1972), 71-103.

\bibitem{HMS} Hajac, P. M., Matthes, R., and Szyma\'nski, W., \emph{Graph $C^*$-algebras and $\zZ_2$-quotients of quantum spheres.} Proceedings of the XXXIV Symposium on Mathematical Physics (Toruń, 2002). Rep. Math. Phys. 51 (2003), no. 2-3, 215-224.

\bibitem{HMS2} Hajac, P.M., Matthes, R., and Szyma\'nski, W., \emph{Quantum real projective space, disc and spheres},
Algebr. Represent. Theory 6 (2003), 169-192.

\bibitem{HS} Hong, J.H., and Szyma\'nski, W., \emph{Quantum spheres and projective spaces as graph algebras},
Comm. Math. Phys. 232 (2002), 157-188.

\bibitem{HJLM} Han, D., Jing, W., Larson, D., and Mohapatra, R., \emph{Riesz bases and their dual modular frames in Hilbert $C^*$-modules}. J. Math. Anal. Appl. 343 (2008), no. 1, 246-256. 

\bibitem{Hat} Hatcher, A., Algebraic Topology, Cambridge University Press, New York, 2001.

\bibitem{Ji} Ji, R., \emph{On Crossed Product $C^*$-Algebras Associated with Furstenberg Transformations on Tori}, PhD Thesis, State University of New York, Stony Brook, 1986.

\bibitem{KR1} Kadison, R.V., and Ringrose, J.R., Fundamentals of the theory of operator algebras.
Vol. I. Elementary theory, Graduate Studies in Mathematics, 15. Amer. Math. Soc., Providence, RI, 1997.

\bibitem{KR2} Kadison, R.V., and Ringrose, J.R., Fundamentals of the theory of operator algebras.
Vol. II. Advanced theory, Graduate Studies in Mathematics, 16. Amer. Math. Soc., Providence, RI, 1997.

\bibitem{KW} Kajiwara, T., and Watatani, Y., \emph{Hilbert $C^*$-bimodules and continuous Cuntz-Krieger algebras}. J. Math. Soc. Japan 54 (2002), no. 1, 35-59.

\bibitem{K} Katsura, T., \emph{On $C^*$-algebras associated with $C^*$-correspondences}. J. of Functional Analysis 217 (2004), 366-401.

\bibitem{K1} Katsura, T., \emph{A class of $C^*$-algebras generalizing both graph algebras and homeomorphism $C^*$-algebras. I. Fundamental results}. Trans. Amer. Math. Soc. 356 (2004), no. 11, 4287-4322.

\bibitem{K2} Katsura, T., \emph{A class of $C^*$-algebras generalizing both graph algebras and homeomorphism $C^*$-algebras. II. Examples}. Internat. J. Math. 17 (2006), no. 7, 791-833.

\bibitem{K3} Katsura, T., \emph{A class of $C^*$-algebras generalizing both graph algebras and homeomorphism $C^*$-algebras. III. Ideal structures}. Ergodic Theory Dynam. Systems 26 (2006), no. 6, 1805-1854.

\bibitem{K4} Katsura, T., \emph{A class of $C^*$-algebras generalizing both graph algebras and homeomorphism $C^*$-algebras. IV. Pure infiniteness}. J. Funct. Anal. 254 (2008), no. 5, 1161-1187.

\bibitem{Kum} Kumjian, A.,\emph{Notes on C*-algebras of graphs,} Contemporary Math. 228 (1998) 189-200.

\bibitem{KP1} Kumjian, A., and Pask, D., \emph{Higher rank graph $C^*$-algebras}, New York J. Math. 6 (2000), 1-20.

\bibitem{KP2} Kumjian, A., and Pask, D., \emph{Actions of $\zZ^k$ associated to higher rank graphs}, Ergodic Theory \& Dynamical Systems, 23 (2003), 1153-1172.

\bibitem{KPR} Kumjian, A., Pask, D., and Raeburn, I., \emph{Cuntz-Krieger algebras of directed graphs}, Pacific J.
Math. 184 (1998), 161-174.

\bibitem{KPRR} Kumjian, A., Pask, D., Raeburn, I., and Renault, J., \emph{Graphs, groupoids, and Cuntz-Krieger
algebras}, J. Funct. Anal. 144 (1997), 505-541.

\bibitem{lan} Lance, E.C., Hilbert $C^*$-modules: A toolkit for operator algebraists, London Mathematical Society
Lecture Note Series, vol. 210, Cambridge University Press, 1995. 

\bibitem{MRS} Mann, M.H., Raeburn, I., and Sutherland, C.E., \emph{Representations of finite groups and Cuntz-Krieger algebras}, Bull. Austral. Math. Soc. 46 (1992), 225–243.

\bibitem{MRS2} Mann, M.H., Raeburn, I., and Sutherland, C.E., \emph{Representations of compact groups, Cuntz-Krieger algebras, and groupoid $C^*$-algebras} in Miniconference on probability and analysis (Sydney, 1991), 135–144, Proc. Centre Math. Appl. Austral. Nat. Univ., 29, Austral. Nat. Univ., Canberra, 1992.

\bibitem{McThesis} McCann, S.J., \emph{$C^*$-algebras associated with topological group quivers}, PhD Thesis, University of Calgary, Calgary, Alberta, Canada, 2012.

\bibitem{Mc1} McCann, S. J., \emph{$C^*$-algebras associated with topological group quivers I: generators, relations and spatial structure}, preprint

\bibitem{MS0} Muhly, P., and Solel, B., \emph{Tensor algebras over $C^*$-correspondences (representations, dilations, and
$C^*$-envelopes)}, J. Funct. Anal. 158 (1998), 389-457.

\bibitem{MS} Muhly, P., and Solel, B., \emph{On the Morita Equivalence of Tensor algebras}, Proc. London Math.
Soc. 81 (2000), 113-168.

\bibitem{MT0} Muhly, P., and Tomforde, M., \emph{Adding tails to $C^*$-correspondences}, Doc. Math. 9 (2004), 79-106.

\bibitem{MT} Muhly, P., and Tomforde, M., \emph{Topological quivers}. Internat. J. Math. 16 (2005), no. 7, 693-755.

\bibitem{M} Munkres, J. R., Topology, 2nd edition, Prentice-Hall, 2000. 

\bibitem{PS} Pask, D., and Sutherland, C.E., \emph{Filtered inclusions of path algebras; a combinatorial approach
to Doplicher-Roberts duality}, J. Operator Theory 31 (1994), 99–121.

\bibitem{P} Paulsen, V.I., Completely bounded maps and operator algebras, Cambridge Studies
in Advanced Math., 78, Cambridge University Press, Cambridge, 2002.

\bibitem{Pims} Pimsner, M., \emph{A class of $C^*$-algebras generating both Cuntz-Krieger algebras and crossed products by $\zZ$}, in Free Probability Theory, fields inst. Commun., vol. 12, Amer. Math. Soc., Providence, 1997, pages 189-212.

\bibitem{Sz} Schweizer, J., \emph{Crossed Product by $C^*$-correspondences and Cuntz-Pimsner Algebras}, `in $C^*$-algebras,
Proceedings of the SFB-Workshop on $C^*$-algebras, Muenster, 1999,' (Eds.) J. Cuntz, S. Echterhoff, Springer Verlag, Berlin, 2000. 

\bibitem{S} Stacey, P.J., \emph{Crossed products of $C^*$-algebras by endomorphisms}, J. Austral. Math. Soc. (Series A) 54 (1993), 204-212.

\bibitem{RW} Raeburn, I., Williams, D.P., Morita Equivalence and Continuous-Trace $C^*$-Algebras, Math. Surveys \& Monographs, vol. 60, Amer. Math. Soc., Providence, 1998.

\bibitem{R} Rieffel, M., \emph{$C^*$-algebras associated with irrational rotations}, Pacific J. Math. 93(2) (1981), 415-429.

\bibitem{RLL} R$\phi$rdam, M., Larsen, F., Lausten, N.J., An Introduction to K-theory for C*-algebras. 256 pp. London Mathematical Society, Student Text 49, Cambridge University Press, Cambridge, 2000. 

\bibitem{RS} R$\phi$rdam, M., St$\phi$rmer, E., Operator Algebras and Non-Commutative Geometry, Vol VII: Classification of Nuclear C*-Algebras. Entropy in Operator Algebras. Encyclopaedia of Mathematical Sciences 126. Springer Verlag, Heidelberg, 2001. 

\bibitem{W} Walters, P., An introduction to ergodic theory, Springer, New York, 1982. 

\bibitem{Yam} Yamashita, S., \emph{Circle Correspondence $C^*$-algebras}, Houston J. Math. 37 (2011), no. 4, 1181-1202.\end{thebibliography}
\end{document}